\input amstex
\documentstyle {amsppt}
\magnification=\magstep1
\hsize=6.5 truein
\vsize=9.0 truein
\NoRunningHeads
\define \OXx{\Cal O_{X,x}}
\define \bk{\bold k}
\define \bC{\bold C}
\define \bR{\bold R}
\define \OX{\Cal O_X}
\define \OtX{\Cal O_{\tilde X}}
\define\op{\operatorname}

\define\cod{\op{cod}}
\define\map#1#2#3{{#1}:(\bk^{#2},0)\rightarrow(\bk^{#3},0)}
\define\MAP#1#2#3{{#1}:{#2}\rightarrow{#3}}
\define\pd#1#2{\frac{\partial#1}{\partial#2}}
\define\tr{\pitchfork}
\define\nottr{\not\tr}
\define \dist{\op{dist}}
\define \em{\frak m}
\define \ord{\op{ord}}
\define\bP{\bold P}
\define\A{\Cal A}
\define\C{\Cal C}
\define\E{\Cal E}
\define\TC{C_T}
\define\PTC{\bold PC_T}

\topmatter
\title Equisingularity of sections, $(t^r)$ condition, and the integral
closure of modules\endtitle
\author Terence Gaffney, David Trotman and Leslie Wilson\endauthor
\affil
Northeastern University, University of Provence,
University of Hawaii
\endaffil
\thanks
\endthanks
\endtopmatter
\document
\par The theory of the integral closure
of modules provides a powerful
tool for studying stratification
conditions which are connected with
limits of linear spaces. It gives an
expression which is both algebraic and geometric
for these conditions.  This connection allows
one to control many geometric
phenomena through associated numerical invariants.
This paper will illustrate these points
by examining the $(t^r)$ conditions which
were introduced by Thom and the second author.
We will show how to apply these
conditions to the study of certain
families of  sections of an analytic space.

The $(t^r)$ conditions deal with the $C^r$ sections of
some stratified set; they
were introduced initially by Thom in 1964 and developed by
the second author from 1976 on,
more recently in collaboration with Kuo and the third author
(\cite{Th}, \cite{Tr1}, \cite{Tr2}, \cite{Tr3},
\cite{Tr4}, \cite{Ku-Tr}, \cite{T-W}), and were applied
to prove various equisingularity results.
For real and complex analytic sets, we show that the
$(t^r)$ conditions have
algebraic formulations in terms of integral closure of modules.
Our formulation gives a new simple
proof, for analytic sets, of the change in the conditions under
Grassmann modification proved by Kuo and the second author \cite{Ku-Tr}
for subanalytic sets; this
is used in conjunction with the principle of
specialization of integral
dependence to give numerical criteria for
familes of plane sections of
complex complete intersections
to be Whitney equisingular.
Some of the results in this paper were announced
by the first author in the 1994
Sao Carlos proceedings \cite{G4}.

In Section 1, we review the notions of integral closure,
reduction, and strict dependence for submodules of $\OXx^p$,
where $X,x$ is the germ of a complex analytic set.  We describe
the analogues which are needed for the case of real analytic sets.
We will apply these concepts to submodules of the Jacobian
module $JM(F)$, where $X=F^{-1}(0)$.  We review results from
\cite{G4} which use these tools to analyse the limits of
tangent hyperplanes to $X,x$ and to characterize the Whitney
$(a)$ and Verdier $(w)$ regularity conditions.

In Section 2, we define jets of transversals and the $(t^r)$
conditions.  The main result of this section is Theorem 2.7,
which is the characterization of condition $(t^r)$ in terms
of integral closure of modules.  We note that condition $(t^0)$
is in fact $(w)$.  The Grassmann modification is a generalization
of the blow-up, in which the projective space is replaced by
the Grassmann space of
planes of the dimension of the transversals we are using.
The Grassmann modification theorem (2.13) says that the modification
improves condition $(t^r)$ to condition $(t^{r-1})$. This is used
at the end of Section 2 and several places in Section 4
for inductive proofs.  Another useful tool is the characterization
in 2.19 of which families of transversals are Verdier equisingular.
This is applied (in 2.28) to families of transversals sharing a common
$r$-jet for which $(t^r)$ holds, and (in 2.23) to
families of transversals sharing a common
$(r-1)$-jet for which $(t^r)$ holds for generic parameter values.
Then we describe the $(t^{r^-})$ condition,
which strengthens $(t^r)$ and was
introduced in \cite{T-W}, and we show that $(t^{r^-})$ holds for
an $(r-1)$-jet iff $(t^r)$ holds for all $r$-jets lying over that
$(r-1)$-jet (assuming $r\ge 1$; $(t^{0^-})$ is the condition
known as strong Verdier regularity). We conclude
by briefly discussing the ambient versions of the $(t^r)$ conditions.

In Section 3, we prove theorems relating conditions $(a)$ and $(t^1)$
to the conormal modification and the aureole of $X$, and use these to
give new examples showing that $(t^1)$ is strictly weaker than $(a)$.

In Section 4, wholly devoted to complex analytic
sets, we state and apply the principle of
specialisation of integral dependence for coherent
submodules of ${\Cal O}^p_X$ where $X$ is
equidimensional, due to  the first author and S.
Kleiman [G-K] (extending the work of Teissier, who
introduced this principle for ideals). We characterise
$(t^r)$ in terms of the genericity of the multiplicity
of a certain submodule of the jacobian module, then
use the principle of specialisation of integral
dependence to give an equimultiplicity criterion for
$(t^r)$. As a consequence we obtain numerical criteria
for Verdier equisingularity of families of plane
sections in various situations.

\head \S 1.  Background on the theory of integral closure of modules
\endhead

In this section we review some basic facts about the
integral closure of a module. We begin with a definition.

\definition{Definition 1.1} Suppose $X,x$ is a
complex analytic set germ, $M$ a submodule of
$\OXx^{p}$. Then $h\in \OXx^{p}$ is in the integral closure
of $M$, denoted $\overline{M}$, iff for all analytic
$\phi :(\bC,0) \rightarrow (X,x)$, $ h \circ \phi
\in (\phi^{*}M){\Cal O}_{1}$.
\enddefinition

An algebraic definition of the integral closure of a
module has been given by David Rees in
\cite{R}, cf. \cite{H-M2}. If we are working in the real analytic case,
we simply use real analytic curves instead
of complex analytic curves.  The set of germs satisfying
the condition of Definition 1.1 in this case
is called the {\it real integral closure} of $M$ and
is sometimes denoted by ${\overline{M}}_{\bR}$, but for simplicity
we will use $\overline{M}$ for real integral closure in this paper.
To see
that these notions are different, it is easy to check
that the integral closure in ${\Cal O}_{{\bC}^2}$ of the
ideal $(x^2+y^2)$ is  $(x^2+y^2)$,
while the integral closure in ${\Cal O}_{{\bR}^2}$
of the ideal $(x^2+y^2)$ is
$(x^2,y^2,xy)$. In the first case, it is clear that $h$
must vanish on the curve defined by $(x^2+y^2)$,
hence must be divisible by $(x^2+y^2)$. In the
second case, it is clear that $\phi^{*}(x^2+y^2)
{\Cal O}_{1}$ is just $(s^{2k})$ where $k$
is the minimum of the orders of the first non-vanishing
term in the components of $\phi$. It is clear that
the same equality holds for $\phi^{*}(x^2,y^2,xy)
{\Cal O}_{1}$.  An algebraic definition of the
real integral closure of an ideal has been
given by Brumfiel \cite{B}. Some of the properties
of the real integral closure can be found in \cite{G1}.
Our results will hold in both the real and complex cases
unless we say otherwise.  Let $\bk$ denote either $\bR$ or $\bC$.

The connection between the integral closures of ideals
and modules is given by
the next proposition, which in fact proves that our
definition agrees with that
of Rees. Let $A$ be the $p\times q$-matrix representing
the left-most homomorphism in any exact sequence
$\OXx^q\rightarrow\OXx^p\rightarrow
{\OXx^{p}}/M\rightarrow 0$.
We denote by $J_k(M)$ the $(p-k)$-th
Fitting ideal of ${\OXx^{p}}/M$, i.e. the
ideal in $\OXx$ generated by the
$k\times k$-minors of $A$; if $h$ is
an element of $\OXx^{p}$,
we denote by $(h,M)$ the module generated by $M$ and $h$.

\proclaim{Proposition 1.2}
Suppose that $M$ is a submodule of
$\OXx^{p}, h\in \OXx^{p}$ and the rank of
$(h,M)$ is $k$ on each component of $X,x$. Then
$h \in \overline{M}$ iff $J_{k}(h,M)
\subseteq \overline{J_{k}(M)}.$
\endproclaim
\demo{Proof} Cf. \cite{1.7} and \cite{1.8}
of \cite{G1}.\qed\enddemo

\comment
As a corollary of 1.2 we obtain:

\proclaim{Corollary 1.3}
Suppose $X,0$ is an equidimensional
complex analytic set defined
by $F :(\bC^N,0) \rightarrow (\bC^p,0)$.
Then $M'\subset JM(F)$ is a reduction of
$JM(F)$ iff for all $\phi : (\bC,0)
\rightarrow (X,0)$ such that
$\phi ({\bC} - \{ 0 \}) \subseteq X - S(X)$, we have
$(\phi^{*}M'){\Cal O}_{1} =(\phi^{*}JM(F)){\Cal O}_{1}$.
\endproclaim
\demo{Proof}
Since $X,0$ is equidimensional (of dimension $d$),
the rank of $JM(F)$
is $N-d$ on each component of $X,0$. The singular
set of $X$, denoted by
$S(X)$, is exactly where the ideal $J_{N-d}(JM(F))$
vanishes. The implication
$\Rightarrow$ is trivial, while the implication
$\Leftarrow$ is now clear,
since if $\phi({\bC}) \subseteq S(X)$
both $J_{N-d}(JM(F))$ and
$J_{N-d}(M')$ vanish.\qed\enddemo
\endcomment

Roughly speaking an element $h$ is in the integral
closure of an ideal $I$  if the
order of vanishing of $h$ on the zero set of $I$  is
the same  as the order of vanishing of $I$. Sometimes
we require the order of vanishing of $h$ to
be greater. The notion of strict
dependence makes this precise.
Suppose $M$ is a submodule of $\OXx^{p},
h \in \OXx^{p}$. Then $h$ is {\it
strictly dependent} on $M$ if for all $\phi : (\bk,0)\rightarrow (X,x)$ we have
$h \circ\phi\in \em_{1}\phi^{*}M$,
where $\em_{1}$ is the
maximal ideal in ${\Cal O}_{1}$. We denote
by $M^{\dagger}$
the set of elements strictly dependent on $M$ .

Note that $M$ in general neither contains nor is
contained in ${{M}}^{\dagger}$. (For
example if $ M = ( x^{3},y^{3} )$ in ${\Cal O}_{2}$,
then ${{M}}^{\dagger} =m^{4}_{2}$.)  If $M$ is an ${\Cal O}_{n}$-module then ${{M}}^{\dagger}$
clearly contains $\em_n{\overline{M}}$ but
this inclusion may be strict.

\comment
By reviewing the proofs of 1.2 and 1.3 we can check
that in the set-up of 1.3,
with $M'$ a reduction of $JM(F)$, then whether or not
$h\in JM(F)$ is strictly
dependent on $M'$ can be checked using curves $\phi$
such that $\phi$ maps
${\bC} - \{ 0 \}$ to the smooth part of $X$.
\endcomment

The precise analogue of 1.2 is:

\proclaim{Corollary 1.3} Suppose $M$ is a submodule
of ${\Cal O}^{p}_{X,x}, h
\in \OXx^{p}$ and the rank of $(h,M)$ is $k$ on each
component of $X,x$. Then
$h \in{{M}}^{\dagger}$ iff each minor in
$J_{k}(h,M)$ which depends on $h$ lies in
${{J_k(M)}}^{\dagger}$.\endproclaim

\proclaim{Proposition 1.4}
Suppose $X,x$ is a real or
complex analytic set germ with irreducible
components $X_i$, $i=1,\dots,l$, $M$ a submodule of
$\OXx^{p}$. In the real case, assume for each $i$ that
the regular points of $X_i$ are dense in $X_i$ in the
metric topology.  Let $U$ be a Zariski open, dense subset of
$X-S(X)$ (hence $U_i=U\cap X_i$ is dense in $X_i$ in the
metric topology.  Then $h\in \OXx^{p}$ is in $\overline{M}$
(respectively, in $M^{\dagger}$),
iff for all analytic $\phi :(\bk,0) \rightarrow (X,x)$
with $\phi(t)\in U$ for all $t\ne 0$, $ h \circ \phi
\in (\phi^{*}M){\Cal O}_{1}$ (respectively, $ h \circ \phi
\in (\phi^{*}M)\em_1$).
\endproclaim
\demo{Proof} We follow the proofs of \cite{1.7}, \cite{1.8} and  \cite{4.2}
of \cite{G1}.  Assume  $h\in \OXx^{p}$ and for all analytic
$\phi :(\bk,0) \rightarrow (X,x)$
with $h(t)\in U$ for all $t\ne 0$, $ h \circ \phi
\in (\phi^{*}M){\Cal O}_{1}$.  Let $S=X-U$ and $S_i=X_i-U$. For the case of integral closure, it suffices to show that
if $\phi$ is a curve in $S$, then $ h \circ \phi
\in (\phi^{*}M){\Cal O}_{1}$. Let $\phi$ be such a curve;
necessarily $\phi$ lies in some $S_i$. 
If $X_i$ is singular, let $(\tilde X_i,\pi)$ be a smooth resolution
of $X_i$; then we may as well assume $\tilde X_i$ is
the germ at 0 of $\bk^m$, $m=\dim X_i$.
Since the regular points of $X_i$ are dense in $X_i$,
$\pi(\tilde X_i)=X_i$.  Hence $\phi$ lifts to a curve $\tilde \phi$ to $\tilde X_i$.

For any $k$ there is a curve $\tilde\phi_1$ so that
$\tilde\phi-\tilde\phi_1\in \em_1^{k+1}{\Cal O}_{1}^p$, and so that
$\phi_1(t)=\pi\circ\tilde\phi_1(t)\in U$ for all $t\ne 0$. Note that $\phi-\phi_1\in \em_1^{k+1}{\Cal O}_{1}^p$.

Suppose $\phi^*(M)\subsetneq
\phi^*((h,M))$.  Then by Artin-Rees and Nakayama's Lemma,
there exists $\nu_o$ such that for all $k>\nu_0$
$$
\phi^*(M){\Cal O}_{1}\subsetneq\phi^*((h,M)){\Cal O}_{1}
\text{ mod } \em_1^k{\Cal O}_{1}^p.
$$

Assume $k>\nu_0$.
Then $\phi_1^*((h,M))=\phi_1^*(M)$ by assumption, hence
$\phi^*((h,M))=\phi^*(M)$ mod $\em_1^k{\Cal O}_{1}^p$,
which is a contradiction.

For the case of strict integral closure, the same argument
works, replacing $\phi^*((h,M))$ by $(\phi^*h,\em_1\phi^*M)$.
\qed\enddemo

In the applications we make, we are interested
in an  
analytic space $X = F^{-1}(0)$ defined by
$F :(\bk^N,0) \rightarrow (\bk^p,0)$.
The modules we use are the jacobian module,
denoted by $JM(F)$, which is the submodule
of $\OXx^{p}$ generated by the first order partial
derivatives of $F$, various submodules of it,
and their reductions. (A
submodule $M'$ of $M$ is a {\it reduction}
of $M$ if $\overline{M'} = \overline{M}$.)

To see the connection between integral closure
and geometry, suppose that the function
$\MAP F{\bC^{n+1},0}{\bC ,0}$
defines the germ of a hypersurface $X$ at
$0$ in $\bC^{n+1}$. Suppose $X$
contains the $y$-axis, $\{0\}\times\bC$.
When does every sequence of
limiting tangent planes to $X$ at the
origin contain the $y$-axis?
In \cite{G3} it is shown that for each limiting
tangent plane $H$ to a hypersurface $X$
one can find a curve $\phi(t)$,
such that the image of $\phi$ lies in the smooth set
of $X$ except for $\phi(0)$, and the
limit as $t$ tends to 0 of the tangent
planes to $X$ at $\phi(t)$ is H.
This means that
$$
\mathop{\lim}\limits_{t\rightarrow 0}
(1/t^{k}) (DF(\phi(t)) =
(a_{1},\ldots,a_{n+1})
$$
where $ a_{1}x_1+\ldots +a_nx_n+a_{n+1}t=0$ defines
$H$, and $k$ is the minimum of
the orders of the first non-vanishing terms in
$\pd F{x_1}\circ \phi,\ldots,\pd F{x_n}\circ
\phi,\pd F t\circ \phi$.  The
condition that $H$ contains
the $y$-axis is just that $a_{n+1}=0$.
This is equivalent to asking that
$\pd Ft\circ \phi\in \em_1\phi^*((\pd F{x_1},
\ldots,\pd F{x_n})\Cal O_{n+1})$.
This means that the condition that every limit of
tangent planes to $X$ at the origin
contain the $y$-axis is equivalent to $\pd Ft\in  ((\pd
F{x_1},\ldots,\pd F{x_n})\Cal O_{n+1})^{\dagger}$.

We can also describe the
limiting tangent planes in terms of integral
closure.  If $H$ is a limiting
tangent plane, and
$v=(v_1,\ldots,v_{n+1})$ is any vector in
$H$, then $\sum a_iv_i=0$, so the
order in $t$ of $DF(\phi(t))(v)$
must be greater than $k$.  Let $JM(F)_{H}$
denote the submodule of $JM(F)$ generated by
$\frac{\partial F}{\partial v}$
for all $v\in H$. Then we see that $H$ is a
limiting tangent plane iff $JM(F)_{H}$ is
not a reduction of $JM(F)$. Both
of these results hold more generally than
for hypersurfaces, with an important change. In the hypersurface case, the Jacobian matrix of
$F$ had only one row, and this
determined the tangent plane to $X$ at a smooth
point of $X$. In general, the rows of the
Jacobian matrix of $F$ will each
determine a tangent hyperplane
to $X$ at a smooth point, and the row space of the
Jacobian matrix of $F$ at a
smooth point will determine
all of the tangent hyperplanes. Thus our
results in the general case will be
stated in terms of tangent hyperplanes.
Here are the corresponding results.

\proclaim{Proposition 1.5} (Theorem 2.2 of \cite{G3})
Suppose $X$ is an equidimensional complex
analytic set, defined by a map germ F.
A hyperplane $H$ is a limiting tangent
hyperplane to $X$ at $0$ iff
$JM(F)_{H}$ is not a reduction of $JM(F)$.
\endproclaim

\proclaim{Proposition 1.6} (Theorem 2.3 of \cite{G3})
Suppose $F$ and $X$ are as above.
Then every limiting tangent hyperplane
contains a given vector
$v$ iff $DF(v)$ is in ${JM(F)}^{\dagger}$.
\endproclaim

Integral closure can also be used to
describe the tangent cones of curves.  Given
$\MAP \phi{\bC}{\bC^n}$, then the tangent cone to
the image of $\phi$ at the origin will be a line, and
the next proposition gives a
test for a plane to contain this line.

\proclaim{Proposition 1.7} Suppose $\phi$ is as above,
and $P$ is a plane. Then the
tangent line to $\phi$ at the origin
lies in $P$ iff $\phi^*(I(P))\subseteq \em_1\phi^*\em_n$.
\endproclaim
\demo{Proof} Let $k$ be the minimum of the
orders of the first non-vanishing term in
the components of $\phi$. Then
$\phi^*\em_n=(t^k)$. We have that
$$
\mathop{\lim}\limits_{t\rightarrow 0} (1/t^{k})
(\phi(t)) = (v_{1},\ldots,v_n).
$$
If $h=a_{1}x_1+\ldots +a_nx_n$ is a generator
of $I(P)$, then the tangent direction
$(v_{1},\ldots,v_n)$ at 0 lies in $P$ iff
$$
\mathop{\lim}\limits_{s\rightarrow 0}
(1/t^{k})(h(\phi(t)))=
a_1v_{1}+\ldots+a_nv_n=0
$$
iff the order of vanishing of every element of
$\phi^*(I(P))$ is greater than $k$.
\qed\enddemo

Two important conditions in the study of
stratifications are Whitney's
conditions $(a)$ and $(b)$. Whitney's condition
$(a)$ holds for a pair of strata
($X$,$Y$) at a point $y$ if every limit of tangent
hyperplanes to $X$ at $y$ contains
the tangent space to $Y$  at $y$. This condition
seems to be necessary for any reasonable condition of
equisingularity along a stratum $Y$ (as opposed to the weaker
$(t^r)$ conditions which give equisingularity of families of sections
through a fixed point $y$). For
example, if
$\{X_s\}$ is any topologically
trivial family of complex  analytic hypersurfaces
with isolated singularities, and $X$ is the total
space of the family, with $S$
the singular set, then ($X-S$, $S$) satisfies
Whitney $(a)$ \cite{L-S}. It is an open question,
posed in particular by Thom, whether topologically
trivial families of complex analytic sets are
always Whitney $(a)$-regular over the smooth
parameter space.

Whitney's condition $(b)$ holds at $y$ if
for any sequence of pairs of points
$(x_i,y_i)$ which converge to ($y$, $y$),
such that the secant line joining the
pair converges to a line $l$, and the tangent
plane to $X$ at $x_i$ converges
to a plane $T$, then $T$ contains $l$. This
condition implies local topological triviality
but it also preserves some important local
infinitesimal structure, the
aureole \cite{L-T}.  For further discussion
of these conditions see \cite{GWPL} and \cite{T3}.
We are interested in a third regularity
condition, Verdier's condition $(w)$ \cite{V}.
Roughly speaking, this condition says that as you approach $y$ from $X$, the distance
between the tangent space to $X$ at $x$ and
the tangent space to $Y$  at $y$ goes to zero
at least as fast as the distance between $x$ and $Y$
goes to zero. In the complex analytic case
$(w)$ has been proved to be equivalent to
Whitney's condition $(b)$ (\cite{T3}, \cite{H-M1}),
while in the real analytic case it implies $(b)$,
but there are real algebraic examples found
by Brodersen and the second author \cite{B-T} which show
that $(b)$ can hold even though $(w)$ fails.
(See \cite{Tr3} for a discussion of the
implications between the various conditions.)  We are
interested in $(w)$ because integral closure methods
connect more directly with $(w)$  than with $(b)$
in the real analytic case. The precise definition is as
follows:

\definition{Definition 1.8}  Suppose A, B are linear
subspaces at the origin in ${\bC}^{N} $, then let
$$
\dist(A,B) = \mathop {\sup}\limits_{\matrix u\in
{B^{\perp}}-\{0\}\\ {v\in {A-\{0\}}}\endmatrix}
\tan \arcsin{ {|(u,v)|} \over {\left\| u \right\|
\left\| v \right\|}}.
$$
\enddefinition

In the applications $A$ is the ``small'' space and $B$ the ``big'' space. Note that $\dist(A,B)=0$ iff
$A\subseteq B$ and $\dist(A,B)=\infty$ iff $A$ has
a nonzero vector perpendicular to $B$.
Note also that if $A$ and $B$ are lines meeting
in an angle $\theta$, then $\dist(A,B)=\tan \theta$.
This distance allows us to talk about the
Whitney condition $(a)$ holding with
a certain exponent.

\definition{Definition 1.9}  Suppose
$y \in Y \cap \bar X$, where $X$
and $Y$ are strata in a
stratification of an analytic space such that
$ \dist (T_yY, T_xX) \le C \dist (x,Y)^e$ for some
constant $C > 0$.  Then
$(X,Y)$ is said to satisfy Whitney $(a)$
with exponent $e$ at $y \in Y$.
Verdier's condition $(w)$ is  Whitney
$(a)$ with exponent 1.
\enddefinition

Verdier proved in \cite{V} that $(X,Y)$ is
locally topologically trivial if $(w)$ holds and
$X\cup Y$ is locally closed.

The usual definition of $\dist(A,B)$ omits the
term $\tan \arcsin $; the definition above was
introduced in \cite{T-W}; it is needed to make
(2.4) below, which represents condition $(t^r)$,
work correctly in the case $r=0$, where we want
(2.4) to represent condition $(w)$.

 In this paper we let $X_0$
denote the smooth points of $X$ (if $X$ is an
analytic space given as $X=F^{-1}(0)$, then
we mean that $X$ is smooth  and
the component functions of $F$
define the reduced structure at points of $X_0$). The
theory of the integral closure of modules
allows us to show:

\proclaim{Theorem 1.10} Suppose $X,0
\subseteq {\bC}^{N}$ is an equidimensional complex
analytic set, $X = F^{-1}(0)$, $Y$ a smooth
subset of $X$. Then Whitney $(a)$
holds for the pair $(X_{0},Y)$ at the origin iff
$\frac{\partial F}{\partial y}$ is strictly
dependent on $JM(F)$ for all tangent vectors
$\frac{\partial}{\partial y}$ to $Y$ at the origin.
\endproclaim
\demo{Proof} Cf. \cite{G3} Cor. 2.4.
\qed\enddemo

\proclaim{Theorem 1.11}  Let $X,Y$ be as above
with coordinates chosen so that
${0} \times {\bC}^k = Y$, $\em_n=(x_1, \dots,x_n)$
denoting the ideal defining $Y$, and let
$\MAP F{{\bC}^N}{\bC}^p$ define $X$ with reduced
structure. Then ${{\partial F} \over
{\partial y}}\in \overline {\em_{n}({{\partial F}
\over {\partial x_1}},\dots{{\partial F}
\over {\partial x_{n}}})\OX}$  for all tangent
vectors ${\partial } \over {\partial y}$ to
$Y$ iff $(w)$ holds for $(X_0,Y)$.
\endproclaim
\demo{Proof} Cf. \cite{G1} Theorem 2.5.
\qed\enddemo

The analogous results hold in the real analytic
case using the real integral closure instead.

In Definition 1.8, suppose that $B$ is a
hyperplane, the kernel of
$\omega \in \op{Hom}(\bC^N,\bC)$.  Then
$$
\dist(A,B) = \mathop {\sup}\limits_{v\in {A-\{0\}}}
\tan \arcsin{{|\omega(v)|} \over {\left\|
\omega \right\|\left\| v \right\|}}.
$$
If $A(t)$ and $B(t)=\ker\omega(t)$ are
analytic in $t$, and $v_1(t),\dots,v_N(t)$ are
analytic and form a basis of $\bC^N$ for all
$t$ small such that $v_1(t),\dots,v_n(t)$ is a
basis of $A(t)$ for all $t$ small, then the order of vanishing
$$
\ord(\dist(A(t),B(t)))=\ord\left(\dfrac
{\max_{1\le i\le n} \omega(t)(v_i)}
{\max_{1\le j\le N} \omega(t)(v_j)}\right),
$$
which is non-negative and, if positive,
$$
\ord(\dist(A(t),B(t)))=\min_{1\le i\le n}
\ord(\omega(t)(v_i))
-\min_{n+1\le j\le N}\ord( \omega(t)(v_j)). \tag 1.12
$$
We will use this formula in the proof of
the main theorem of the next section.

\head \S 2. The $(t^r)$ conditions and the Grassmann
modification
\endhead

\definition{Definition 2.1} Suppose $X$ and
$Y$ are disjoint smooth submanifolds
of ${\bR}^N$ in a neighborhood of 0 $\in Y$. For $r$ a positive integer,
we say $X$ is $(t^r)$ regular over $Y$ at 0
when every $C^r$ submanifold $Z$
transverse to $Y$ at 0, and of complementary
dimension to $Y$, is transverse to $X$ near 0.
\enddefinition

This concept was first introduced by Thom \cite{Th};
Thom allowed $Z$ to be of dimension
greater than the complementary dimension of $Y$.
Transversals $Z$ of dimension
complementary to $Y$ are called
{\it direct transversals}.  The second author realized that
the degree of smoothness of the transversal was
important: $(t^r)$ implies $(t^{r+1})$, but
the converse does not hold in general.

Further
the second author showed that if you work with direct
transversals, although Whitney $(a)$ implies $(t^1)$,
the converse is only known to hold when $X$ and $Y$
are subanalytic and the dimension of $Y$ is 1 (\cite{Tr1},
\cite{Tr2}, \cite{Tr3}, \cite{Tr4}). In
\cite{Tr4} the second author gave an example of a semialgebraic
$X$ and $Y$, with $\dim Y = 2$,
satisfying $(t^1)$ but not $(a)$.
We will give algebraic examples of this phenomenon in
Section 3, where we relate conditions $(a)$ and $(t^1)$
to the conormal modification and the aureole of an analytic
set $X$.

A little thought reveals the connection between these
ideas and the family of
sections of $X$ by planes
which are direct transversals to $Y$. If, as we move
through the family of slices,
one of the planes has a higher order of contact
with $X$ than the others, than that slice should not be
equisingular with the
other slices. The failure of $(t^1)$ for
a plane is a way of detecting this higher contact.

It was shown in \cite{T-W} that it is convenient
to refine the concept of $(t^r)$ to one in which
we only require that those transversals to $Y$ which
share a common $r$-jet be transverse to $X$; we
give the appropriate definition of $r$-jet of a
transversal next.

\definition{Definition 2.2} First suppose $r>0$. If $\MAP{f,g}{(\bold R^n,0)}{(\bold R^k,0)}$ are
$C^r$ map-germs with $j^rf(0)=j^rg(0)$, we say
that $P=\Gamma(f)$ and $Q=\Gamma(g)$
are $r$-equivalent. The equivalence class
is called the $r$-jet of $P$, denoted $j^rP$.
Note that this condition is equivalent to
$|f(x)-g(x)|=o(|x|^r)$. Now suppose $r=0$.
We want our transversals to have well-defined tangent
planes away from the origin, so we don't want to consider
graphs of continuous functions when forming a 0-jet.
Instead we look at graphs of mapgerms $f$ of class $C^{0,-1}$, that is,
$f$ is $C^0$,   $C^1$ except possibly at 0, and
$|f'(x)|=O(|x|^{-1})$.  For two such $f$ and $g$, we will say
that $j^0f(0)=j^0g(0)$ if $f(0)=g(0)$
and $|f'-g'|=o(|x|^{-1})$; as above we say
that $P=\Gamma(f)$ and $Q=\Gamma(g)$
are $0$-equivalent.
The graph of the
$r$-th Taylor polynomial of $f$ is called the
degree $r$ polynomial representative of $j^rP$, denoted $P_r$.
If $\MAP{f}{(\bold C^n,0)}{(\bold C^k,0)}$, then the
$r$-jet of $P=\Gamma(f)$ is the equivalence class of $f$ as
a real map $\MAP{f}{(\bold R^{2n},0)}{(\bold R^{2k},0)}$.
\enddefinition

\remark{Remark} Assume that $f$ and $g$ are $C^1$
map-germs, and $|f(x)-g(x)|=o(|x|^r)$, $r>0$. Then
$\dist(z,\Gamma(g))=o(|z|^r)$ for all $z\in\Gamma(f)$.
Conversely, if $\dist(z,\Gamma(g))=o(|z|^r)$
for all $z\in\Gamma(f)$, then there exists a
$\MAP{\phi}{(\bold R^n,0)}{(\bold R^n,0)}$ such that
$$
|(x,f(x))-(\phi(x),g(\phi(x)))|=o(|(x,f(x))|^r)
\text{ and hence }=o(|x|^r),
$$
since the derivative of $f$ is bounded near 0. Hence,
$$
|x-\phi(x)|\text{ and } |f(x)-g(\phi(x))|=o(|x|^r).
$$
But since the derivative of $g$ is bounded near 0,
$$
|g(x)-g(\phi(x))|\le C|x-\phi(x)|=o(|x|^r)
\text{ for some }C>0.
$$
Thus, $|f(x)-g(x)|=o(|x|^r)$.

One can define $r$-equivalence of germs $P$ and
$Q$ of $n$-manifolds by requiring that
both $\dist(z,Q)=o(|z|^r)$ for all $z\in P$
and $\dist(z,P)=o(|z|^r)$ for all $z\in Q$
hold.  This notion of $r$-equivalence is clearly
a $C^1$ invariant. So $r$-equivalence
as in (2.2) is the restriction of this notion
to graphs of $C^r$ map-germs.
\endremark

\definition{Definition 2.3} Suppose $\MAP{f}{(\bold R^n,0)}{(\bold R^k,0)}$
is $C^r$, $r> 0$, or $C^{0,-1}$ (in the latter
case let $r=0$), $P=\Gamma(f)$, and $W$ is a submanifold of ${\bold R}^{n+k}$
containing 0 in its closure. Then $W$ is $(t^r)$ for
$P$ if for all $Q$ such that $j^rQ=j^rP$ (i.e. they
are $r$-equivalent) then $Q\tr W$ near (but not
necessarily at) $0$.
\enddefinition

If $W$ is a complex analytic submanifold of $\bold C^{n+k}$, then we can take
$f: (\bold R^{2n},0)\to  (\bold R^{2k},0)$, so $P\subset \bold R^{2n+2k}$
and the above definition applies to $W$ as before.

Later on we will give another approach to the $(t^r)$ conditions in the analytic case.


\remark{Remark} Suppose $ Y=\{0\}\times \bk^k$ and $P$ is the graph of the
constant function $f(x)=0$. It is shown in \cite{T-W} that the
pair $(W,Y)$ satisfies
Verdier's condition $(w)$ at $(0,0)$ iff  $W$ is $(t^0)$  for $P$.
In this case, we will also say
that $W$ is $(t^0)$ regular over $Y$ at $0$.
\endremark
\medpagebreak

By \cite{T-W}, for $r\ge 0$, a manifold $W$ disjoint from $Y=\{0\}\times \bk^k$  fails to be $(t^r)$ for
$P=\Gamma(f)$ if and only if 
\medpagebreak

\noindent(asfd):  there exist $a_i=(x_i,y_i)\in W
\subseteq \bk^{n+k}$, such that
$x_i\rightarrow 0$ and there exist $n$-dimensional
planes $T_i\nottr T_{a_i}W$ such that
$$
\align |y_i-f(x_i)| &= o (|x_i|^r)\tag {2.4.1}\\
\dist(T_{b_i}P,T_i) &= o (|x_i|^{r-1}),\tag {2.4.2}
\endalign
$$
where $b_i=(x_i,f(x_i))$.

Clearly $T\nottr T_aW$ iff there exists a
hyperplane $H\supseteq T+T_aW$. If
$T\subseteq H$, $\dist(A,T)\ge \dist(A,H)$. Given $A$ and $H$, there exists
$T\subseteq H$ with $\dim T=\dim A$ such that
$\dist(A,T)=\dist(A,H)$.  From
this, it follows that (2.4) is equivalent to:\newline
\smallpagebreak

$x_i\rightarrow 0$ and
hyperplanes $H_i\supseteq T_{a_i}W$ such that
$$
\align |y_i-f(x_i)| &= o (|x_i|^r)\\ 
\dist(T_{b_i}P,H_i) &= o (|x_i|^{r-1}) 
\endalign
$$
where $b_i=(x_i,f(x_i))$.

If, in addition, $W$ is subanalytic and $P$ is analytic, an
application of the Curve Selection Lemma
(applied as in Theorem 5.3 of \cite{T-W}) shows
that (2.4) is equivalent to:\newline
\smallpagebreak

\noindent (2.5):
there exists an analytic curve $\phi(t)=(x(t),y(t))$ such
that $\phi(t)\in W$ for $t\ne 0$ and $\phi(0)=0$,
and there exists an analytic curve of
hyperplanes $H(t)\supseteq T_{\phi(t)}W$ such that
$$
\align |y(t)-f(x(t))| &= o (|x(t)|^r)\tag {2.5.1}\\
\dist(T_{b(t)}P,H(t)) &= o (|x(t)|^{r-1})\tag {2.5.2}
\endalign
$$
where $b(t)=(x(t),f(x(t)))$.

\definition{Definition 2.6}   Let $Y=\{0\}\times \bk^k$.  Suppose
 $P=\Gamma(f)$ is a $C^r$ (or $C^{0,-1}$) direct transversal to $Y$,
$r> 0$ (or $r=0$).
Suppose $X$ is an analytic subset of $\bk^{n+k}$
containing 0 with singular set $S(X)$, and
$S$ is a closed subset such that $S(X)\cup (Y\cap X)\subseteq S\subseteq X$.
Let $W=X-S$; note that $W\cap Y=\emptyset$. Then we say $(X,S)$ is $(t^r)$ for $P$ if
$W$ is $(t^r)$ for $P$ and 
$$
\text{for all $Q$ such that $j^rQ=j^rP$, then $Q\cap S=\{0\}$}
\tag{2.6.1}$$
(Often the $S$ will be understood and we will simply
say that $X$ is $(t^r)$ for $P$, or even $(t^r)$ holds for $P$; $S$ is assumed to be
$S(X)\cup (X\cap Y)$ if we don't say otherwise.)
\enddefinition

\comment
\remark{Remark} Suppose $y\in Y$ and $P$ is the graph of the
constant function $f(x)=y$. Then $X$ is $(t^0)$  for $P$ at $y$ iff
the germ of $S$ at $y$ is contained in $Y$, and the pair $(W,Y)$ satisfies
Verdier's condition $(w)$ at $(0,y)$. In this case, we will also say
that $X$ is $(t^0)$ (or $(w)$) regular over $Y$ at $y$.
\endremark
\endcomment
\remark{Remark} Suppose that $P$ is the graph of the
constant function $f(x)=0$. Then $X$ is $(t^0)$  for $P$ at $0$ iff
the germ of $S$ at $0$ is contained in $Y$, and the pair $(W,Y)$ satisfies
Verdier's condition $(w)$ at $(0,0)$. In this case, we will also say
that $X$ is $(t^0)$ (or $(w)$) regular over $Y$ at $0$.
\endremark
\medpagebreak
Here is an example on $\bR^4$ with coordinates $(x_1,x_2,y_1,y_2)$.
Let $F(x_1,x_2,y_1,y_2)= (x_1- y_1^2-y_2^2)^3-x_2^2$, let $X=F^{-1}(0)$
and let $W=X_0$. Then $W$ is $(t^0)$ regular over $Y$ at 0, but
$X$ is not $(t^0)$ regular over $Y$ at 0 ($S$ is the surface
$x_1= y_1^2+y_2^2$, $x_2=0$).  Note $W$ being $(w)$ regular over $Y$ at 0
does not imply $Y\subseteq X$, although it can be shown that
$X$ being $(w)$ regular over $Y$ at 0
does  imply $Y\subseteq X$. 
In some of our applications we will be interested in
proving results about the family of sets $\{X\cap P_a\}$
where $P_a$ is a family of direct transversals, so we
want to control the intersection of $P_a$ with all of $X$,
not just $X-S$.  One way to accomplish this is to stratify
$X$ ($W$ being the top stratum) such that each stratum is
$(t^r)$ for $P$. We discuss this briefly after Cor 2.29.
A   special  case of this is to require that
$(X,S)$ be $(t^r)$ for $P$ (as defined above); then the transversals will miss
all the strata in $S$ (except at 0) and hence satisfy the
$(t^r)$ condition.  In the complex case, this is exactly
the case where the dimension of $S$ is less than or
equal to the codimension of $P$. In the real case, we can
have $S$ of larger dimension than that but still have
transversals $P$ missing $S$ away from 0.

It is obvious that $X$ being $(t^r)$ for $P$
depends only on the $r$-jet of $P$. Thus we may
talk about $X$ being $(t^r)$ for $z$, where $z$
is an $r$-jet of some $P$.  Further
observe that disjoint manifolds $W$ and $Y$ are $(t^r)$ regular over $Y$ at 0
(as in Definition 2.1)   if,
and only if, $X=W\cup Y$ is $(t^r)$ for all $C^r$ direct
transversals to $Y$ through 0; in the complex case we say
$X$ is {\it $(t^r)$ regular over $Y$ at 0} if $X$ is $(t^r)$
for all $r$-jets of graphs of complex polynomials vanishing
at 0. (In this complex analytic setting, the curve of 2.5 can be
taken to be complex analytic.) We say $X$ is $(t^r)$ regular over $Y$ at $y$ if
$\{x-y \, |\, x\in X\}$ is {\it $(t^r)$ regular over $Y$ at 0}.

If $\bk=\bC$, then $X$ can be $(t^r)$ over $Y$
in the real or in the complex
sense: the distinction is whether one requires
$(t^r)$ over all $r$-jets, or just
over $r$-jets of graphs of complex polynomial
mappings from $\bk^n$ to $\bk^k$ of degree $\le r$. However, we do not know of any example of an
analytic $X\cup Y$ so that $X$ is $(t^r)$ over
$Y$ in the complex but not in the real sense.

The next theorem provides our algebraic criterion
for condition $(t^r)$. To motivate this criterion,
consider analytic map-germs $f$ and $g$, with
$f(0)=0$ and the components of $g$ in $m_n^r$. Let $f_t=f+tg$; the family of transversals
$P_t=\Gamma(f_t)$ are $(r-1)$-equivalent.  By
results of \cite{T-W}, if $X=F^{-1}(0)$
and $Y=\{0\}\times\bk^k$ contains the points of $X$ at
which $F$ is singular, and if $(X,Y\cap X)$ is $(t^r)$
for $P=\Gamma(f)$, then there is a family
of homeomorphism-germs $h_t:P\rightarrow P_t$,
$h_t(0)=0$, preserving $X$, for $t$ sufficiently
small.  Temporarily we ask instead for diffeomorphism
germs, and require them to preserve all fibers of $F$. In fact let us suppose somewhat more: assume that the
$h_t$ are diffeomorphism germs on $(\bk^{n+k},0)$,
preserve $Y$ and the fibers of $F$, and $h_t(P)=P_t$
for all $t$.  Let $\tilde P=\cup_tP_t\times\{t\}$
for all $t$ and let $H(x,y,t)=(h_t(x,y),t)$.
Then $H$ is the flow of the germ of a
vector field $\xi_t(x,y)+\pd{}{t}$ on $k^{n+k+1}$
which is tangent to $\tilde P$ and to $Y\times \bk$
such that $DF(\xi_t)=0$ for all $t$.

Let $\Theta(\bold k^{n+k},P,Y)$ denote the
$\Cal O_{n+k}$-module
of germs of vector fields which are tangent
to $P$ and to $Y$. It is not hard to show
that the requirements placed on
$\xi$ imply that $\xi_0=g \mod  \Theta(\bold k^{n+k},P,Y)$.
Let $I(P)$ denote the ideal of analytic function-germs
vanishing on $P$ ($I(P)=
\{y_1-f_1,\dots,y_k-f_k\}\Cal O_{n+k}$ if $P=\Gamma(f)$).
We let $JM(F)_P$ denote the submodule of $JM(F)$
generated by $\pd{F}{x_i}+
\sum_j\pd{f_j}{x_i}\pd{F}{y_j}$, $i=1,\dots,n$;
this submodule is
obtained by applying to $F$ vector fields tangent to the
fibers of the submersion $y-f(x)$ defining $P$.
(Note that if $P$ is a plane this module is the
module we denoted by $JM(F)_P$ previously.)
If $(x_1,\dots,x_n,y_1,\dots,y_k)$ are coordinates
on ${\bk}^{n+k}$,
let $JM_y(F)$ denote the submodule of $JM(F)$ generated
by the partials of $F$ with respect to $y$.  Hence for this
analytic equisingularity  to hold for all $g$, we need that
$$
\em_n^rJM_y(F)\subset
DF(\Theta(\bold k^{n+k},P,Y))=m_nJM(F)_P+I(P)JM_y(F)$$
(equality of these two modules is an easy exercise).
 If one only requires that
the diffeomorphisms preserve $X$, then one replaces
these $\Cal O_{n+k}$-modules by
the corresponding $\Cal O_X$-modules. To formulate our algebraic criterion
for $(t^r)$, we replace the module on the
right-hand side of this inclusion
by its integral closure, as in the next Theorem.

\medskip
\definition{Definition 2.6 (continued)}
We say that a smooth submanifold $W\subset X\subset\bk^{n+k}$ is defined by   $\map F{n+k}p$ and $Y=0\times\,\bk^{k}$ if $X=F^{-1}(0) $ and
$W$ is the Zariski open subset  of $X-Y$ at which
$F$ has rank $c= \cod (X-Y)$ (possibly less then $p$,
to allow for non-complete intersections).
\enddefinition

In the complex analytic  case, if the Zariski closure of
$X-Y$ is equidimensional, and $F$ has  generic rank $
c= \cod (X-Y)$ on each component of the Zariski closure of $X-Y$, then the set of points $W$ at which $F$ has
rank $c= \cod  (X-Y)$ is a metric dense subset of $X-Y$
(in fact of $X$). The real  analytic case is different. Consider
$F(x,z,y)=z^2-y^2x$. The zero set of $F$, $X$ is the
Whitney umbrella  with handle the $x$ axis. Let $Y$ be the
$y$ axis, so the metric  closure of $X-Y$ is $X$. Then $F$
and $Y$ do not define a metric  dense $W$, as the metric
closure doesn't contain the handle. (Of  course, if we
choose $Y$ to be the handle, then $W$ is metric dense  in $X-Y$.)

In relating the ($t^r$) conditions to the Grassmann modification
in the  real analytic case, we will assume $W$ is metric dense
in $X-Y$; and  this condition will be preserved by Grassmann
modification, whereas  the condition that the smooth points of
$X$ be metric dense in  general won't.

Note that $W$ consists of regular points, but we don't insist
that it be all regular points. Let $\Cal O_{X,0}$ be the
ring defined by quotienting by the ideal generated by
the components of $F$.  Let $S=X-W$; so this
contains all singular points of $X$ as well as regular
points at which $F$ doesn't define the reduced structure,
as well as all points of $X\cap Y$.  (By construction, $W$
contains  no points of $Y$, hence no points of $X\cap Y$.)
We will usually  assume that $W$ is metric dense in $X-Y$. The germs of $X$, $F$,
etc\. induce germs of sheaves.

\proclaim{Theorem 2.7} Suppose $X\subseteq {\bk}^{n+k}$
is the germ of an analytic space at the origin
defined by $F$, and $Y$, $S$ and $W$ are as in the paragraph above,
$W$ metric dense in $X-Y$, $P$ an analytic direct
transversal to $Y$. For the $(t^0)$ case, assume in addition that
$Y \subset X$.  For $r\ge 0$, $W$ is $(t^r)$ for the transversal $P$  iff
$$
\em_n^rJM_y(F)\subseteq \overline{
{\em_n}JM(F)_P+I(P)JM_y(F)} \tag 2.7.1
$$
(take the integral closure inside $\Cal O_{X,0}^p$,
where $p$ is the number of components of $F$;
in the real analytic situation we take real
integral closure).

Furthermore, $X$ is $(t^r)$ for the transversal $P$ iff
$W$ is $(t^r)$ for the transversal $P$ and
$$
\em_n^r\Cal O_{S,0}\subseteq \overline{I(P)\Cal O_{S,0} }. \tag 2.7.2
$$
\endproclaim

First we prove several preliminary results.  We always
assume $r\ge 0$ unless we state otherwise.

\proclaim{Lemma 2.8}  A curve $\phi(t)=(x(t),y(t))$ satisfies
(2.5.1) iff $\phi^*I(P)\subseteq \em_1 \phi^*\em_n^r$.
\endproclaim
\demo{Proof} Let $\ord(u(t))$ denote the order of the
power series $u(t)$; if $u(t)$ is
vector valued, use the minimum order of the components,
equivalently the order of
$|u(t)|$; if $I$ is an ideal, $\ord(I)$ denotes the
minimum order of all members of $I$.

Suppose $\ord(x(t))=l$.  Then $\ord(|x(t)|^r)=lr=
\ord(\phi^*\em_n^r)$ and $\ord(|y(t)-f(x(t))|)=
\ord(\phi^*I(P))$.  The result is immediate.
\qed\enddemo

\proclaim{Lemma 2.9} Let $\phi(t)$ be a curve satisfying (2.5.1),
with $\phi(t)$ in $W$ when $t\neq 0$.
Then the following are equivalent:
$$
\align
\phi^*(\em_n^{r-1}JM_y(F))\subseteq \phi^*(JM(F)_P)
&\text{ (if $r\ge 1$) or }\\
\phi^*(JM_y(F))\subseteq \phi^*(\em_nJM(F)_P)
&\text{ (if $r=0$)} \tag 2.9.1
\endalign
$$
and
$$
\phi^*(\em_n^{r}JM_y(F))\subseteq \phi^*(\em_nJM(F)_P
+ I(P)JM_y(F)).\tag 2.9.2
$$
\endproclaim
\demo{Proof}
Clearly (2.9.1) implies (2.9.2).  Assume (2.9.2)
holds.  By Lemma (2.8),
$$
\phi^*(I(P)JM_y(F))\subseteq \em_1\phi^*(\em_n^{r}JM_y(F)).
$$
Therefore,
$$
\align  \phi^*(\em_nJM(F)_P+I(P)JM_y(F))&\subseteq\\
\phi^*(\em_nJM(F)_P)+ \em_1\phi^*(\em_n^{r}JM_y(F))&\subseteq\\
\phi^*(\em_nJM(F)_P)+ \em_1\phi^*(\em_n JM(F)_P&+I(P)JM_y(F))
\text{ by (2.9.2).}
\endalign
$$
By Nakayama's Lemma,
$$
\phi^*(\em_n^r JM_y(F))\subseteq
\phi^*(\em_n\phi^*JM(F)_P). \tag 2.9.3
$$
If $r=0$, we are done. Assume $r\ge 1$.
For $l=\ord(\phi(t))$, $\phi^*(\em_n^{r}JM_y(F))=
t^{rl}\phi^*(JM_y(F))$ and $\phi^*(\em_nJM(F)_P)=
t^l\phi^*(JM(F)_P)$.  Thus (2.9.1) holds. \qed
\enddemo

\proclaim{Proposition 2.10} Condition (2.9.1) is
equivalent to requiring that
there does not exist a curve of hyperplanes
$H(t)\supseteq T_{\phi(t)}W$ ($t\ne 0$) satisfying
(2.5.2).
\endproclaim
\demo{Proof} Assume $r>0$;
a slight modification of the argument works for $r=0$.

By Proposition 1.11 of \cite{G1} and the remark after
Proposition 4.2 of \cite{G1},
condition (2.9.1) is equivalent to
$$
\phi^*(\em_n^{r-1}\psi JM_y(F))\subseteq
\phi^* \psi JM(F)_P
\text{ for all analytic } \psi(t)\in
\op{Hom}(\bk^p,\bk).
$$
Letting $v_i(t)=\pd{}{x_i}+\sum_j
\pd{f_j(x(t))}{x_i} \pd{}{y_j}$,
$i=1,\dots,n$ and $w_j=\pd{}{y_j}$,
$j=1,\dots,k$, the above is equivalent to
$$
\forall \psi(t)\in \op{Hom}(\bk^p,\bk),\
\min_{j=1,\dots,k}\ord(\phi^*\em_n^{r-1}
\psi DF\circ \phi \cdot w_j)\ge
\min_{i=1,\dots,n}\ord(\psi DF\circ \phi \cdot v_i).
$$
This fails to hold iff
$$
\exists \psi(t)\in \op{Hom}(\bk^p,\bk),\
\min_{j=1,\dots,k}\ord(\phi^*\em_n^{r-1}
\psi DF\circ \phi \cdot w_j)<
\min_{i=1,\dots,n}\ord(\psi DF\circ \phi \cdot v_i).
$$
If this latter happens, then $\psi DF\circ \phi $ is not identically 0.
Then $H(t):= \ker(\psi DF \circ \phi)$ is a hyperplane
containing $T_{\phi(t)}W$ and,
applying 1.12 to $\omega(t)=\psi DF\circ \phi$, we have
that the above inequality is equivalent to
$$
\dist(T_{b(t)}P,H(t))=  o(|x(t)|^{r-1}),
$$
since $b(t) = (x(t),f(x(t)))$ and $T_{b(t)}P$ is spanned by
$\{v_1(t), \dots ,v_n(t)\}$.
\qed
\enddemo

\demo{Proof of the Theorem} Recall that the failure of $W=X-S$  to
be  $(t^r)$ for $P$  is
equivalent to (2.5) holding for some curve (call it $\phi$) in $W$.

If (2.5.1) fails for   curve $\phi$, then $\phi^*\em_n^r
\subseteq\phi^*I(P)$, which
implies that (2.9.2) holds for  $\phi$.  If (2.5.1) holds
but (2.5.2) fails,
then Proposition 2.10 and Lemma 2.9 imply
that (2.9.2) holds for that $\phi$.

Thus $W$  being $(t^r)$ for $P$ implies that (2.5) does not hold for any analytic curve $\phi$ in $W,0$, which implies that (2.9.2) holds for all analytic curves $\phi$ in $W,0$.   If $\phi$ is
a curve in $Y\cap X,0$ and $r>0$, then (2.9.2)
holds trivially.  If we let $r=0$, then (2.9.2)
holds for such a curve since  $Y\subseteq X$ implies $\phi^*(JM_y(F)) =0$.
Hence by Proposition 1.4 the
integral closure condition of the
Theorem holds.

Now assume (2.5) holds. Then Proposition 2.10 and
Lemma 2.9 imply that (2.9.2) fails to hold
for the curve $\phi$ of (2.5).  Thus the integral
closure condition (2.7.1) fails to hold.

The proof of the statement in the second paragraph
of (2.7) is much easier.
Recall from Definition 2.6 that $X$ $(t^r)$ over $P$ means $W$   $(t^r)$ over $P$ and the   $r$-jet of  $P$  misses $X-W$ except at 0.
The condition (2.7.2)
means:  given any curve $\phi(t)=(x(t),y(t))$ in $S$,
there exists a $C>0$ such that $|x(t)|^r\le C|y(t)-f(x(t))|$
for small enough $t$.   This inequality implies that
the same inequality holds  for any $g$ such that $j^rg(0)=j^rf(0)$,
and hence the $r$-jet of $P$ misses $S$ except at 0 (if $r=0$, the
inequality implies that $S\subseteq Y$ near 0). The converse follows
from the Interpolation Lemma (2.6) of \cite{T-W}. 
\qed\enddemo

If $X$ is $t^1$ for $P$, then this implies that no
line in the tangent cone to $S$ is a tangent line to $P$.

In case $r=0$, letting $P=\bk^n\times\{0\}$, Theorem 2.7 says
that $(w)$ is equivalent to
$$
\pd F{y_i}\in \overline{(\em_n JM_x(F)+\em_kJM_y(F))},
\qquad i=1,\dots,k.
$$
By Nakayama's Lemma, this is equivalent to
$$
\pd F{y_i}\in \overline{\em_n JM_x(F)},
\qquad i=1,\dots,k.
$$
Thus Theorem 1.11 is the special case of Theorem
2.7 in which $r=0$.

In  the real analytic, $(t^0)$ case it is possible
for the integral  closure condition of the theorem
to hold even if $Y\cap X$ is a  point. This is true
of the example that appears after definition 2.6. There, if $\phi(t)$ is a curve on $X$, then the
order in $t$ of $ \phi^*J(F)_y$ is greater than or
equal to the order of $\phi^*(x_1x_2) \in \phi^*m_YJ(F)_z$.
It seems likely that this cannot happen in the  complex analytic case.

It is important to note that Definition 2.6
has two arguments---in applications we can
vary both $P$ and $X$ (or $F$).
This is why the same condition can be used to
study the stratification of a
complex analytic set (take $F$ to be
the defining equations for $X$, $P$
a ``probing set"), or the order of
$V$-determinacy of a map-germ
$\MAP f{{\bC}^n}{\bC}^p$.
To do this, work on the ambient
space of the graph of $f$,
take $P$ to be the graph of $f$, and let
$F=(y_1,\dots,y_p)$; $X$ then becomes
${\bC}^n$, while $Y$ is still given by
$(x_1=\dots=x_n=0)$. If we restrict
$I(P)$ to $X$ we get
$f^*(\em_p)\OX$, while $JM_y(F)$ becomes
$\Cal O_n^p$. The condition of
Definition 2.3, via Theorem 2.7, then becomes
$$
\em_n^r\Cal O_n^p\subseteq {\overline
{\em_nJM(f)+f^*(\em_p)\Cal O_n^p}}.
$$
In \cite{G1}, Corollary 4.6, it is shown that
if this condition holds for $f$,
then any deformation of $f$ which
fixes the $r$-jet of f gives a Verdier equisingular
deformation of $f^{-1}(0)$, hence $f$ is $r$ $V$-determined.

Next we will discuss the Grassmann modification of $X$, and
how it affects the $(t^r)$ condition.

\definition{Definition 2.11} Let $G_{n+k,n}$ denote the
Grassmannian of $n$-planes containing the origin in
${\bk}^{n+k}$, $E_{n+k,n}$ the canonical bundle
associated to  $G_{n+k,n}$; the  fiber  of
$E_{n+k,n}$ over a point $P$
is just the set of points of ${\bk}^{n+k}$
in $P$. Denote the projection
of $E_{n+k,n}$ to ${\bk}^{n+k}$ by $\beta_{n+k,n}$.
If $X$ is a subset of ${\bk}^{n+k}$, we
call $\tilde X =\beta_{n+k,n}^{-1}(X)$, the
$G_{n+k,n}$ modification of $                                                                 X$. When $n$ and $k$
are clear from the context, we will simply refer to the
Grassmann modification of $X$.
\enddefinition

This notion was introduced in connection
with the $(t^r)$ conditions by Kuo
and the second author \cite{Ku-Tr}.
Note that  $G_{n+k,n}$ is embedded in
$E_{n+k,n}$ as the  zero section of
$E_{n+k,n}$. This means that,
if $X$ is an analytic set and $0\in X$, we
can think of $0\times G_{n+k,n}$ as a
stratum of $\tilde X$; note that the
projection to $0\times G_{n+k,n}$
makes $\tilde X$ a family of analytic
sets with $0\times G_{n+k,n}$ as the
parameter space. The members of
this family are just
$\{X\cap P \}$ as $P$ varies through the points
of $G_{n+k,n}$. This means that
if we want to study the equisingularity
of the family $\{X\cap P \}$, then we should
study the regularity conditions
that the pair $({\tilde X}_0,0\times G_{n+k,n})$
satisfy.  In this paper,
we will be interested only in the behavior
of this pair at points of $G_{n+k,n}$
which correspond to direct transversals;
we denote this set by $\tilde Y$.
Kuo and the second author proved a remarkable result
in \cite{Ku-Tr} about the behavior of
regularity conditions under Grassmann
modification. In the following result, $(t^0)=(w)$.

\proclaim{Theorem 2.12}(Kuo-Trotman) Suppose $X$ and $Y$ are disjoint smooth submanifolds of
${\bR}^N$ in a neighborhood of the origin, $Y$ is analytic, and $X$ is
subanalytic. Then, for each $r\ge 0$,
X is $(t^{r+1})$ regular over $Y$  at 0 iff $\tilde X$
is  $(t^r)$ regular over $\tilde Y$ at every point
of $\tilde Y$.
\endproclaim
\demo{Proof} Cf. Theorem 1 in \cite{Ku-Tr} and \cite{T-W}
(there a more general statement is proved,
without the hypotheses that $Y$ be analytic and
$X$ subanalytic).
\qed\enddemo

This result shows in particular that if
$X$ is $(t^1)$ regular  over $Y$
at the origin,  $X\cup Y$ is locally closed,
and $Y$ contains the singular set of $X\cup Y$,
then the family $\{X\cap P \}$ is
Verdier equisingular (i.e. $(w)$ holds),
hence is topologically trivial.

We will prove a version of Theorem 2.12
that holds in both the real and complex
analytic cases using our integral closure
criterion.  In the Grassmann modification of
${\bk}^{n+k}$ we will always be working in
a neighborhood of a direct transversal $P$.
(Recall  we regard $P$ as a point of the zero
section of $E_{n+k,n}$.) Since we are always
working in $E_{n+k,n}$, we denote the projection
$\beta_{n+k,n}$ by $\beta$.  Since all planes
near $P$ are also direct transversals,  they are
also graphs over the $x$-plane, and they
have equations $\{y_1=\mathop{\sum}\limits_j
a_{1,j} x_j,\dots,y_k=\mathop{\sum}\limits_j
a_{k,j} x_j\}$.
This means we have local coordinates on
$E_{n+k,n}$ given by $(x_1,\dots,x_n,a_{1,1},
\dots,a_{k,n})$. In these coordinates we have
$$
\beta (x_1,\dots,x_n,a_{1,1},\dots,a_{k,n})=
(x_1,\dots,x_n,\mathop{\sum}\limits_j a_{1,j}
x_j,\dots,\mathop{\sum}\limits_j a_{k,j} x_j).
$$

Given $\phi :({\bk},0) \rightarrow (\tilde X,\{0\}\times P)$,
then $\beta\circ\phi$
is tangent to $P$ at the origin.

\proclaim{Theorem 2.13}  Suppose $X\subseteq {\bk}^{n+k}$
is the germ of an analytic space at the origin
defined by $F$, and $Y$, $S$ and $W$ are
as in the paragraph before Theorem 2.7. Then,
for each $r\ge 0$, $W$ (respectively $X$) is $(t^{r+1})$
regular over $Y$ at 0 iff $\tilde W$
(respectively $\tilde X$) is $(t^r)$ regular over $\tilde Y$ at every
point of $\tilde Y$.  Moreover, every $(r+1)$-jet $z$ of
direct transversals to $Y$ at 0 lifts to a non-unique
$r$-jet $\tilde z$ of a direct transversal to $\tilde Y$,
and $W$ (respectively $X$) is $(t^{r+1})$
regular for $z$ iff $\tilde W$ (respectively $\tilde X$) is $(t^r)$  for $\tilde z$.
\endproclaim
\demo{Proof}  First note that, since $0\in X$,
$\tilde Y=\beta^{-1}(0)\subseteq \tilde X$.

Since $F$ defines $X$, $G:=F\circ \beta$
defines $\tilde X$. From the chain rule we note that
$$
\pd G{a_{i,j}}=x_j\pd F{y_i}\circ\beta,\text{ and }
\pd G{x_j}=\pd F{x_j}\circ\beta+
\mathop{\sum}\limits_i a_{i,j}\pd
F{y_i}\circ\beta.
$$
Let $\MAP p{\bk^n}{\bk^{nk}}$ be a polynomial
map of degree $\le r$; $\Gamma(p)$
is a direct transversal to $\tilde Y$. Then $\beta(\Gamma(p))=\Gamma(q)$
for a polynomial mapping $q$ of degree
$\le r+1$, $q(0)=0$,
$$
q=(\sum_j p_{1,j}x_j,\dots,\sum_j p_{k,j}x_j).
$$
Conversely every such $q$ has $\Gamma(q)=
\beta(\Gamma(p))$ for some
(not unique) $p$ of degree $\le r$.

We need to show that $\tilde W$ is $(t^r)$
for $\Gamma(p)$ iff $W$ is $(t^{r+1})$ for
$\Gamma(q)$.  Let $\MAP \phi{(\bk,0)}{(\tilde W,p(0))}$
be an analytic curve which has order $r+1$
contact with $\Gamma(p)$:
$$
\phi(t)=(x(t),a(t)),\quad |a(t)-p(x(t))|=o(|x(t)|^r).
$$
Then $\phi_1=\beta\circ \phi$ has order
$r+2$ contact with $\Gamma(q)$:
$$
\phi_1(t)=(x(t),y(t)),\quad y(t)=(\sum_j a_{1,j}(t)x_j(t),
\dots,\sum_j a_{k,j}(t)x_j(t)),
$$
so
$$
y(t)-q(x(t))=(\sum_j (a_{1,j}(t)-p_{1,j}(x(t)))x_j(t),
\dots)=o(|x(t)|^{r+1}).
$$

By the proof of Theorem 2.7, we know that the inclusion
$$
{\phi_1}^*(\em_n^{r+1}(JM_y(F)))\subseteq
{\phi_1}^*(\em_n(JM(F)_{\Gamma(q)})
+I(\Gamma(q))(JM_y(F)))\tag 2.14
$$
(which is the module condition defining $(t^{r+1})$
along the curve $\phi_1$) is equivalent to
$$
{\phi_1}^*(\em_n^rJM_y(F))\subseteq {\phi_1}^*
(JM(F)_{\Gamma(q)}).\tag 2.15
$$
The vector fields tangent to $\Gamma(q)$ are generated by
$\pd{}{x_j}+\sum_i\pd{q_i}{x_j}\pd{}{t_i}$,
$j=1,\dots,n$, so (2.15) is  equivalent to
$$\multline
{\phi}^*(\beta^*\em_n^r\pd F{y_i}\circ \beta)\subseteq
{\phi}^*(\pd F{x_j}\circ\beta+
\sum_i\pd{q_i}{x_j}\pd F{y_i}\circ\beta)\Cal O_1=\\
{\phi}^*(\pd F{x_j}\circ\beta+
\sum_i(\sum_l\pd{p_{i,l}}{x_j}x_l+p_{i,j})
\pd F{y_i}\circ\beta)\Cal O_1
\endmultline\tag 2.16
$$
(and note that $\beta^*\em_n^r=\em_n^r$).

Similarly the module condition for $(t^r)$
along $\phi$ reduces to
$$
\phi^*(\em_n^{r-1}JM_a(G))\subseteq\phi^*
(JM(G)_{\Gamma(p)}).\tag 2.17
$$

Now $JM(G)_{\Gamma(p)}$ is generated by
$$
(\pd {}{x_j}+
\sum_{i,l}\pd{p_{i,l}}{x_j}\pd{}{a_{i,l}})(F\circ\beta)=
\pd F{x_j}\circ\beta+\sum_i a_{i,j}\pd F{y_i}\circ\beta+
\sum_{i,l}\pd{p_{i,l}}{x_j}x_l\pd F{y_i}\circ\beta
$$
and
$JM_a(G)$ is generated by $\pd G{a_{i,j}}=
x_j\pd F{y_i}\circ\beta$.

Thus the left hand sides of (2.16) and (2.17)
are identical.  The difference
of the right hand sides of (2.16) and (2.17)
is generated by
$$
\phi^*(\sum_i(a_{i,j}-p_{i,j})\pd F{y_i}\circ\beta)=
\sum_i(a_{i,j}(t)-p_{i,j}(x(t))) \phi^*\beta^*\pd
F{y_i}\in \em_1\phi^*(\em_n^r\beta^*JM_y(F))
$$
for $j=1,\dots,n$.  So if either inclusion (2.16) or
(2.17) holds, by Nakayama's Lemma both must
hold and the two modules are in fact the same. So if we assume $(t^{r+1})$ for $\Gamma(q)$,
then (2.16) holds for all curves $\phi_1$ which
have order $r+2$ contact with $\Gamma(q)$,
which implies (2.17) holds for
all curves $\phi$ having order $r+1$ contact with
$\Gamma(p)$, and hence we have
$(t^r)$ for $\Gamma(p)$.

The converse is similar.

The proof that $X$ is ($t^{r+1}$) iff $\tilde X$ is ($t^r$)
is similar, but much simpler, and we omit it.\qed\enddemo

It should be noted that, even if $F$ defines the reduced
structure on $X$, the Grassmann modification $G$ of $F$
may not define the reduced structure on $\tilde X$.
For example, suppose $X$ is the curve $F(x,y,z)=(y,x^3+z^2)=0$.
Let $\beta(x,y,u,v)=(x,y,ux+vy)$.  Then the ideal generated
by $F\circ \beta$ is $(y,x^2(x+u^2))$, which gives a
non-reduced structure on the component $\tilde Y$ of $\tilde X$.

The Grassmann modification will be a major tool
in section 4.  Next we will treat
an important tool for studying equisingularity
of families of transversals.

Suppose $X\subseteq {\bk}^{n+k}$
is the germ of an analytic space at the origin
defined by $F$, and $Y$, $S$ and $W$ are
as in the paragraph before Theorem 2.7. Let $\map {f(x,u)}{n+c}k$
determine a germ of a family of direct transversals
$\Gamma(f_u)$, where $f_u(x)=f(x,u)$. Assume $f(0,u)=0$ for all $u$. Let
$\beta(x,u)=(x,f(x,u))$, $G=F\circ \beta$ and
$\tilde X=\beta^{-1}X=G^{-1}(0)$; let $\tilde S=\beta^{-1}S $ and
$\tilde W=\beta^{-1}W $; and $\tilde Y=\beta^{-1}Y$
is the germ at (0,0) of $\{0\}\times \bk^c$. ( Note that $\tilde Y\subseteq \tilde X$; note
that the germ of the map to the base at a point of $\tilde Y$ of the
Grassmann modification is a special case
of such a $\beta$).
\definition{Definition 2.18}
The family $\Gamma(f_u)$ has {\it Verdier
equisingular intersection with $W$} (respectively, {\it with $X$})
if $\tilde W$ is $(w)$ regular over
$\tilde Y$ at each point of $\tilde Y$
(respectively, $\tilde X$ is $(w)$ regular
over $\tilde Y$ at each point of $\tilde Y$).
\enddefinition

Recall that $\tilde X$ being $(w)$ regular over $\tilde Y$
at each point of $\tilde Y$ implies, in particular, that
$ \tilde S-\tilde Y$ misses a neighborhood
of $\tilde Y$ (this follows from (2.6.1));
since $\tilde W\cup\tilde Y$ is then locally closed in
this neighborhood of $\tilde Y$, it further implies that there is a family of rugose (with respect to 0)
homeomorphisms $h_t$ from $\Gamma(f_0)$
to $\Gamma(f_u)$ preserving $X$.

\proclaim{Theorem 2.19 (of Thom-Levine type)}
The family $\Gamma(f_u)$ has Verdier
equisingular intersection with $W$  iff,
for all $i=1,\dots,c$,
$$
\left(\pd Fy\circ \beta\right)\cdot \pd f{u_i}\in
\overline{\em_n (dF\circ\beta)_* JM_x(\beta)},
\tag 2.20
$$
where $\pd Fy\circ \beta$ takes its values in
$\op{Hom}(\bk^k,\bk^p)$, and $(dF\circ\beta)_*
JM_x(\beta)$ is the $\Cal O_{\tilde X}$-module
generated by composing $dF\circ\beta$ with
elements of $JM_x(\beta)$.

In addition,  $\Gamma(f_u)$ has Verdier
equisingular intersection with $X$  iff
it has Verdier
equisingular intersection with $W$
and $\beta^{-1}S- \tilde Y$ misses a neighborhood of $\tilde Y$.
\endproclaim
\demo{Proof} 
Pick any analytic curve $\phi(t)=(x(t),u(t))$
in $\tilde X$ such that $\phi(0)=0$.
Then (2.20) says that, for all $i=1,\dots,c$
and for all such $\phi$,
$$
\phi^*\left(\pd Fy\circ \beta\cdot \pd f{u_i}\right)\in
\phi^*(\em_n (dF\circ\beta)_* JM_x(\beta) \OtX).
\tag 2.21
$$
The Verdier equisingularity of the intersections of
$\Gamma(f_u)$ with $W$ is equivalent to
$$
\phi^*\pd G{u_i}=\phi^*\left(\pd Fy\circ
\beta\cdot \pd f{u_i}\right)\in
\phi^*(\em_n JM_x(G))=\phi^*(\em_n
(dF\circ\beta)_* JM_x(\beta) \OtX),
\tag 2.22
$$
for all $i=1,\dots,c$.  Clearly (2.21) and (2.22)
are the same.

The last statement of the Theorem follows immediately from the definitions.
\qed\enddemo
\comment
The right hand sides differ by the term
$\phi^*((f-f_0)\OtX JM_y(F)\circ\beta)$.
By the Mean Value Theorem,
$$
\phi^*((f-f_0)\OtX)\subseteq \phi^*(\em_c\OtX JM_u(f))
\in \em_1\phi^*JM_u(f),
$$
so the difference of the right hand sides
is in $\em_1$ times the
module generated by the left hand sides;
by Nakayama's Lemma, (2.22) and
(2.23) are equivalent.
\endcomment


Let $\map {f_0}nk$ be some analytic mapping. Perturb $f_0$ by all homogeneous terms of degree $r>0$:
$$
f_u(x)=f_0(x)+\sum\Sb |\omega|=r\\ 1\le i\le k\endSb
u_\omega^ix^\omega e_i,
$$
where $e_i=(0,\dots,1,\dots,0)$ is the $i$th
standard basis vector in $\bk^k$.
So $f(x,u)=f_u(x): (\bk^n\times\bk^c,0\times\bk^c)
\rightarrow (\bk^k,0)$
is a family of representatives of all $r$-jets
lying over the $(r-1)$-jet determined by $f_0$.
Let $f^i$ denote the $i$-th component function of $f$.

\proclaim{Corollary 2.23} For any $u_0\in\bk^c$,
$W$ (respectively $X$) is $(t^r)$ for $\Gamma(f_{u_0})$ iff the germ
at $u_0$ of the family $\Gamma(f_u)$ has
Verdier equisingular intersection with $W$
(respectively $X$).
\endproclaim
\demo{Proof}  The $(t^r)$ condition for $\Gamma(f_{u_0})$ is:
$$
\em_n^r\OX \pd Fy \subseteq \overline{(\em_n
JM(F)_{\Gamma(f_{u_0})} +
I(\Gamma(f_{u_0}))\pd Fy)\OX} \tag 2.24
$$
(here for convenience we write $\OX \pd Fy$
for $JM_y(F)$).  By Theorem 2.19, Verdier
equisingularity of the intersection of
$\Gamma(f_u)$ with $X$ is equivalent to
$$
\em_n^r \OtX\pd Fy\circ\beta \subseteq
\overline{\em_n \Cal O_{\tilde X}
(dF\circ\beta )\pd{\beta}x}, \tag 2.25
$$
where $\tilde X$ will be taken to mean the
germ of $\tilde X$ at $(0,u_0)$.

Fix a curve $\phi$ in $\tilde X$, $\phi(0)=u_0$,
and let $\phi_1=\beta\circ\phi$ in $X$.
On $\phi_1$, (2.24) becomes
$$
\phi^*(\em_n^r\Cal O_{\tilde X} \pd Fy\circ\beta)
\subseteq \phi^*\left(\em_n \beta^*JM(F)_{\Gamma(f_{u_0})} +
(f^i-f^i_{u_0})\Cal O_{\tilde X}\pd Fy\circ\beta\right), \tag 2.26
$$
and (2.25) is
$$
\phi^*(\em_n^r\Cal O_{\tilde X} \pd Fy\circ\beta)
\subseteq \phi^*\left(\em_n \Cal O_{\tilde X}(\pd Fx\circ\beta
+(\pd Fy\circ\beta)\pd fx)\right). \tag 2.27
$$
Since
$$
(f^i-f^i_{u_0})\circ\phi(t)
=\sum(u_\omega^i\circ\phi(t)-{u_0}_\omega^i)x^
\omega\circ\phi(t)\in \em_1\phi^*\em_n^r,
$$
$$
\left(\pd{f^i-f^i_{u_0}}{x_j}\right)\circ\phi(t)
=\sum(u_\omega^i\circ\phi(t)-{u_0}_\omega^i)
\pd{x^\omega}{x_j}\circ\phi(t)\in \em_1\phi^*\em_n^{r-1},
$$
and since
$$
\beta^*JM(F)_{\Gamma(f_{u_0})}=
\left(\pd Fx\circ\beta+(\pd Fy\circ\beta)
\pd {f_{u_0}}x\right)\Cal O_{\tilde X},
$$
the difference of the right hand
sides of (2.26) and (2.27) is
$$
\phi^*\left(\em_n\Cal
O_{\tilde X} (\pd Fy\circ\beta)\pd{}x(f-f_{u_0}) +
(f^i-f^i_{u_0})\Cal O_{\tilde X}\pd Fy\circ\beta\right)
\in \em_1\phi^*(\em_n^r\Cal O_{\tilde X}
\pd Fy\circ\beta),
$$
which is $\em_1$ times the left hand
side of (2.26) (or (2.27)).
Therefore (2.26) and (2.27) are equivalent.
\qed\enddemo

The next Corollary states that if $X$ is
$(t^r)$ for some $r$-jet, then all
representatives of the jet have Verdier
equivalent intersections with $X$.

\proclaim{Corollary 2.28} Let $\map {h}nk$
be some analytic mapping, and let $g$
be its Taylor polynomial of degree $r$. Let $f_u=(1-u)h+ug$.  If $W$ (respectively $X$) is $(t^r)$
for $\Gamma(h)$,
then $\Gamma(f_u)$ has Verdier equisingular
intersection with $W$ (respectively $X$).
\endproclaim
\demo{Proof}
The proof is a minor variation of that of the
previous theorem.
It suffices to prove the result for the germ
of the family at $u_0$, for each $u_0\in [0,1]$.
 Since $\Gamma(f_{u_0})$ has
the same $r$-jet as $\Gamma(h)$, $X$ is $(t^r)$
for $\Gamma(r_{u_0})$, the proof goes through
with minor notational changes. \qed\enddemo

If $F(x,y)=y$, then $X=\bk^n\times\{0\}$.  Let $S=\{0\}$
and so $W=X-\{0\}$.  Since $\Gamma(f_u)\cap X=f_u^{-1}(0)$,
$\Gamma(f_u)$ having Verdier equisingular
intersection with $X$ is the same as saying `` the
functions $f_u$ are Verdier $V$ equisingular" as defined
in \cite{G1}.

In this case Theorem 2.19 yields
the following result:

\proclaim{Corollary 2.29} (see \cite{G1} for related results).
The family $f_u$ are Verdier $V$ equisingular iff
$$
\pd f{u_i}\in \overline{\em_n JM_x(f)}.
$$
\endproclaim

In this paper we have always assumed that our sections miss
$S(X)$ (or more generally $S$) away from 0.  It is
possible to avoid this restriction.  We will describe this briefly
now, but won't pursue it further in this paper.
Suppose there are mappings $F_i$, $i=0,\dots,l$ defining $X=X^0\supset
X^1\supset \dots \supset X^l$.  Assume $W_i=X^i-X^{i+1}$ is a manifold
containing $X^{i+1}$ in its closure, and the collection of $W_i$'s is
a $(w)$-regular stratification (strata are not necessarily connected).
Then say this stratification $\Cal S$ is $(t^r)$ for $P$ if each $W_i$ is.
We can apply the results we have proved with $F$, $X$, $S$ and $W$
replaced by $F_i$, $X^i$, $S_i=X^{i+1}$ and $W_i$.  The assumption that
$\Cal S$ be $(w)$-regular allows one to apply the results of Section 4
of \cite{T-W} to extend
various theorems of our paper to this stratified case.  For example,
any finite dimension family $Q_t$ of representatives of the $r$-jet of $P$
will have Verdier equisingular intersection with $X$.  The pullbacks
of the $W_i$ by this family will be a family $\tilde W_i$ which is
$(w)$-regular among themselves by the original $(w)$-regularity and its
invariance under transverse intersection, and they are $(w)$-regular over
the $\tilde Y$ by the results of our paper.

A condition $(t^{r^-})$ was introduced in
\cite{T-W}.  It can be defined as for $(t^r)$
using the appropriate notion of $C^{r^-}$
functions, which we will omit. Let $X$, $Y$, $P$, $W$ and
$S$ be as in (2.6).  The failure of $W$ to be
$(t^{r^-})$   for $P$ is equivalent to
(2.4) or (2.5) holding with $o(\text{ })$
replaced by $O(\text{ })$.  If $\map{f,g}nk$
are $C^{r^-}$, so in particular if they are $C^r$,
with $j^{r-1}f(0)=j^{r-1}g(0)$, then
$|f(x)-g(x)|= O(|x|^r)$ and
$\dist(T_x(\Gamma(f)),T_x(\Gamma(g)))=
O(|x|^{r-1})$ (in case $r=0$, one cannot use functions
$f$ which are merely $C^0$---one must also assume they
are $C^1$ except at 0 and that  $\|df(x)\|=
O(|x|^{-1})$; the assumption  $j^{-1}f(0)=j^{-1}g(0)$
is vacuous).
Hence, the analogue of (2.4) implies that $W$
is $(t^{r^-})$ for $\Gamma(f)$
iff it is $(t^{r^-})$ for $\Gamma(g)$.  Furthermore,
$X$ is $(t^{r^-})$ for $P$ iff $W$ is and, for
all $C^{r^-}$ representatives $Q$ of  $j^{r-1}P$, $Q\cap S=\{0\}$ (if we restrict to $C^r$ representatives,
the ``if" of this ``iff" fails). 
Then the analogue of Theorem 2.7 is:

\proclaim{Theorem 2.30} Suppose $X\subseteq {\bk}^{n+k}$
is the germ of an analytic space at the origin
defined by $F$, and $Y$, $S$ and $W$ are as in the paragraph before
(2.7), $W$ metric dense in $X-Y$.
Suppose $P$ is an analytic direct
transversal to $Y$. For the $(t^{0^-})$ case, assume in addition that
$Y \subset X$.  For $r\ge 0$, $W$ is $(t^{r^-})$ for the transversal $P$   iff
$$
\em_n^rJM_y(F)\subseteq
{(\em_nJM(F)_P+I(P)JM_y(F))}^{\dagger} \tag 2.30.1
$$
(take the integral closure inside $\Cal O_{X,0}^p$,
where $p$ is the number of components of $F$;
in the real analytic situation we take real
integral closure).

Furthermore, $X$ is $(t^{r^-})$ for the transversal $P$ iff
$W$ is $(t^r)$ for the transversal $P$ and
$$
\em_n^r\Cal O_{S,0}\subseteq (I(P)\Cal O_{S,0} )^{\dagger}. \tag 2.30.2
$$
\endproclaim

We omit the proof.

\proclaim{Theorem 2.31}  Suppose $X\subseteq {\bk}^{n+k}$
is the germ of an analytic space at the origin
defined by $F$, and $Y$, $S$ and $W$ are
as in the paragraph before Theorem 2.7.  Then,
for each $r\ge 0$, $W$ (respectively $X$) is $(t^{(r+1)^-})$
regular over $Y$ at 0 iff $\tilde W$ (respectively $\tilde X$) is $(t^{r^-})$ regular over $\tilde Y$ at every
point of $\tilde Y$.  Moreover, every $r$-jet $z$ of
direct transversals to $Y$ at 0 lifts to a non-unique
$(r-1)$-jet $\tilde z$ of a direct transversal to $\tilde Y$,
and $W$ (respectively $\tilde X$) is $(t^{(r+1)^-})$
for $z$ iff $\tilde W$ (respectively $\tilde X$) is $(t^{r^-})$ for $\tilde z$
(saying that a condition holds for a $(-1)$-jet
is taken to mean that it holds for all 0-jets).
\endproclaim
\demo{Proof} The proof is a minor variation of
the proof of (2.13). Reduce the degrees of $p$ and $q$
by 1, and replace the Landau symbol $o$ by $O$.
In (2.14)--(2.17), multiply the right hand
sides of the $\subseteq$ by $\em_1$. The rest of
the proof is the same.
\qed\enddemo

\proclaim{Theorem 2.32} If $r\ge 1$,
$W$ (respectively $X$) is $(t^{r^-})$ for $\Gamma(f)$ iff
$W$ (respectively $X$) is $(t^r)$ for all graphs of
degree $\le r$ polynomial representatives
of the ($r-1$)-jet of $f$.
\endproclaim
\demo{Proof} Using the preceding theorem,
we can reduce to the case $r=1$.  This was
proved in \cite{G4}, Corollary 2.15.
\qed\enddemo

\comment
We have already seen that if $X$
is $(t^{r^-})$ for $\Gamma(f)$,
then $X$ is $(t^{r^-})$, and hence $(t^r)$,
for all $C^r$ representatives of the
($r-1$)-jet of $\Gamma(f)$,  This proves
the ``only if\," direction.

Now we prove the ``if\," direction.
Without loss of generality, we can assume $f=0$,
after the invertible polynomial change of coordinates
$$(x,y)\mapsto (x,y-f(x))=(x,y').$$
For simplicity of notation we
continue to use $y$ for $y'$.

Assuming $(t^{r^-})$ fails for $\Gamma(0)$,
then there exist $\phi (t) = (x(t),y(t))$
in $X =  F^{-1}(0)$ and
$\psi(t) \in \op{Hom}({\bk}^p,\bk)$ such that
$$
\phi^*\left(\em_n\psi\pd Fx\OX^n
+\em_k\psi\pd Fy\OX^k\right)
\subseteq \phi^*\left(\em_n^r\psi\pd Fy\right).\tag 2.33
$$
It follows that
$$
\phi^*\em_k
\subseteq \phi^*\em_n^r, \text{ i.e. } y(t)
= O(\vert x(t) \vert ^r), \tag 2.34
$$
and that
$$
\phi^*\psi\pd Fx\OX^n
\subseteq \phi^*\em_n^{r-1}\psi\pd Fy\OX^k.\tag 2.35
$$
If $l=\ord(\phi(t))=\ord(x(t))$, then $\ord(y(t)) \geq lr$ by (2.34).
Also $\ord(\phi^*\em_n^{r-1})=l(r-1)$, so that (2.35)
is equivalent to
$$
\ord(a_i) \geq l(r-1) + \min_{j=1,\dots,k}\ord(b_j)
\text{ for all $i$},
\tag 2.36
$$
where $a_i = \psi \pd{F}{x_i} \circ \phi$ and
$b_j = \psi \pd{F}{y_j} \circ \phi$.

Without loss of generality we can assume that $x(t)$
is tangent to the $x_1$-axis at $0$,
and that $x(t) = (t^l, 0, \dots , 0) + o(t^l)$.
We can also assume that
$$y(t) = (e_1x_1^r, e_2x_1^r, \dots, e_kx_1^r) + o(x_1^r)
= (e_1, \dots , e_k)t^{lr} + o(t^{lr}),$$
where $e_1,\dots,e_k$ are constants, possibly zero.
We then define a homogeneous polynomial $p$ of degree $r$ in
$x_1,\dots,x_n$ by
$$
p(x) = (e_1, \dots , e_k)x_1^r +
x_1^{r-1}L_0\bmatrix x_2\\ \vdots \\x_n\endbmatrix.
$$

To specify $L_0 \in M_{k,{n-1}}(\bk)$ we first write
$$
L^*(t) = \frac{-1}{\sum_{m=1}^k b_m^2}
\bmatrix a_2b_1 &\cdots &a_nb_1 \\
\vdots &   &\vdots \\
a_2b_k &\cdots &a_nb_k \endbmatrix \in M_{k,{n-1}}(\bk\{t\}).
$$

Observe that the entries of $L^*(t)$ all have order at least
$l(r-1)$ in $t$.
Let $L(x_1) = L^*(x_1^{1/l})$.
Then $L(x_1) = O(\vert x_1 \vert ^{r-1})$.
Finally we set $L_0 = (j^{r-1}(L))(1),$ where the
$(r-1)$-jet is taken with respect to the variable $x_1$.
This means that each entry of the $k \times (n-1)$ matrix
$L_0$ is the coefficient of $(t^{l(r-1)})$ in the corresponding
entry of the matrix $L^*(t)$.

We observe that
$$
p(x(t)) = (e_1, \dots , e_k)x_1^r(t) +
L^*(t)\bmatrix x_2\\ \vdots \\x_n\endbmatrix +
\sum_{q=2}^nx_q(o(t^{l(r-1)})).\tag 2.37
$$

To prove that $(t^r)$ fails to hold, we need to show that,
for each $i$ such that $1 \leq i \leq n$,
$$
\ord\left(a_i + \sum_{j=1}^k \pd{p_j}{x_i}b_j\right)
 > l(r-1) + \min_{j=1,\dots,k}\ord(b_j).\tag 2.38
$$
First suppose $i \geq 2$.

By (2.33),
$$
p_j(x(t)) = e_jx_1^r -
\frac{b_j}{\sum_{m=1}^kb_m^2}\sum_{q=2}^na_qx_q +
\sum_{q=2}^nx_q(o(t^{l(r-1)})).
$$
Hence
$$
\pd{p_j}{x_i}(x(t)) = - \frac{a_ib_j}{\sum_{m=1}^kb_m^2}
+ o(t^{l(r-1)}).
$$
It follows that if $s = \min_{j=1,\dots,k}\ord(b_j)$, then
$$
a_i + \sum_{j=1}^k \pd{p_j}{x_i}b_j =
a_i - \frac {a_i\sum_{j=1}^k b_j^2}{\sum_{j=1}^k b_j^2}
+ o(t^{l(r-1)s}) = o(t^{l(r-1)s}),
$$
proving (2.38)
for $i \geq 2$.

Now we take $i = 1$.

Because $\phi(t) \in F^{-1}(0)$,
$\phi'(t)$ is tangent to
$F^{-1}(0)$ at $\phi(t) = (x(t),y(t))$, i.e.
for all $\psi(t) \in \op{Hom}(\bk^p,\bk)$,
$\phi'(t) = (x'(t),y'(t))$ is orthogonal to
$(\psi \pd{F}{x_1},\dots,\psi \pd{F}{x_n},\psi
\pd{F}{y_1},\dots,\psi \pd{F}{y_k})\circ \phi(t)$.

In particular,
$$
\sum_{i=1}^n x_i'a_i + \sum_{j=1}^k y_j'b_j = 0.
$$
Hence
$$
a_1 = - \frac{(\sum_{i=2}^n x_i'a_i +
\sum_{j=1}^k y_j'b_j)}{x_1'}.
$$
Then
$$
a_1 + \sum_{j=1}^k \pd{p_j}{x_1}b_j =
- \sum_{i=2}^n \frac{x_i'}{x_1'}a_i
+ \sum_{j=1}^k \left(\pd{p_j}{x_1}
- \frac{y_j'}{x_1'}\right)b_j.
$$
Because $\ord(\frac{x_i'}{x_1'}) =
\ord(\frac{x_i}{x_1}) > 0$ if $i \geq 2$,
the first term on the right has order in $t$
strictly greater than $\ord(a_i)$.

So by (2.36), the first term on the right has order
$> l(r-1) + \min_{j=1,\dots,k}\ord(b_j)$,
as required in (2.38).

Now it suffices to prove that for each $j$,
$1 \leq j \leq k$,
$$
\ord\left(\pd{p_j}{x_1} - \frac{y_j'}
{x_1'}\right) > l(r-1) = \ord(x_1^{r-1}).\tag 2.39
$$
By definition of $p$ and the expression for $y(t)$,
$$
\align
\left(\pd{p_j}{x_1} - \frac{y_j'}{x_1'}\right) &=
\left(re_jx_1^{r-1} + (r-1)x_1^{r-2}L_0
\bmatrix x_2\\ \vdots
\\x_n\endbmatrix \right)
 - \left(re_jx_1^{r-1}\frac{x_1'}{x_1'} +
o(t^{l(r-1)})\right) \\
&= (r-1)x_1^{r-2}L_0\bmatrix x_2\\ \vdots \\
x_n\endbmatrix + o(t^{l(r-1)}),
\endalign
$$
which has order $ > \ord(x_1^{r-1}) = l(r-1)$ because
$\ord(x_i) > \ord(x_1)$ if $i > 1$.

This proves (2.39) and hence
(2.38), and completes the proof that the failure of
$(t^{r^-})$ for $\Gamma(f)$ implies that
$(t^r)$ fails for some
$C^r$ representative of the $(r-1)$-jet of $f$.
\qed\enddemo
\endcomment


Condition $(t^{0^-})$ turns out to be the so-called ``strong Verdier regularity", where
$X$ is strong Verdier over $Y$ at $0$ if
for every $\epsilon > 0$ there is an ambient neighbourhood of $0$ on which
$\dist(T_xX,T_yY) < \epsilon \dist(x, Y)$. This condition was called
``differentiably regular" in \cite{K-T-L}. It is a non-generic condition and
implies a rather strong
local topological triviality---the trivialising
homeomorphisms, $C^1$ off $Y$,
have a differential on $Y$ which is
the identity in the normal direction to $Y$  (shown in \cite{T-W}), so that for
example the density of the singular variety $X$ is {\it constant} along $Y$
\cite{Va1}. Verdier equisingularity merely ensures that the density
is lipschitz along $Y$, as shown by G. Valette \cite{Va2}
(G. Comte \cite{C} previously proved continuity).

\bigskip

There is also the notion of the ambient or relative $t^r$
conditions.

{\it Set-up for ambient $(t^r)$.}  Let $X\subset{\bk^{n+k}}$ be an
analytic space defined by a function $G$, and
let $f$ be a function mapping $X$ to $\bk$. Pick an extension of $f$ to
the ambient space, which we also denote by $f$.
Denote the mapping with components $(f,G)$ by $F$. Let $X(f)$ denote
$F^{-1}(0)$. Let $S=S(f,G)\cup \{0\}\times\bk^k$.

In what follows the results are independent of the chosen extension of
$f$.

We use the same definition of the $r$-jet of a $C^r$ manifold as before.

\definition{Definition 2.33} In the ambient $(t^r)$ set up,  $f$ is
$(t^r)$ for
$P$ if for all $Q$ such that $j^rQ=j^rP$ (i.e. they
are $r$-equivalent) then $Q$ misses $S$ near (but not at) 0
and is transverse to the fibers of
$f:X-S\rightarrow \bk$ near (but not
necessarily at) $0$.
\enddefinition

It is possible to define a notion of ambient $(t^r)$ for $f:X\rightarrow
\bold {k}^p$, but in the complex case it rarely holds.

To see this consider a map-germ $f:\bold {C}^3\rightarrow \bold {C}^2$,
whose singular set $S(f)$ is a non-empty curve. Let $P$
be the germ of any smooth complex analytic surface at the origin which
intersects $S(f)$ only at
the origin. Then the restriction of $f$ to $P$ is not a submersion, hence
$S(f|P)$ is at least a curve. At each point of
$S(f|P)$ different from the origin, $P$ is tangent to the fiber of $f$ at
that point. So  $(t^r)$ fails for such
$P$ for all  $r$.

If $Y$ is a submanifold passing through the origin,
then we say that $f$ is $(t^r)$ regular over $Y$ at $0$ if
and only if $f$ is $(t^r)$ for all $C^r$ direct
transversals to $Y$ through $0$.
In the complex case we say that
$f$ is $(t^r)$ regular over $Y$ at $0$ if $f$ is $(t^r)$
for all $r$-jets of graphs of complex polynomials vanishing
at 0.

Let $X_f$ denote the subset of $X$ consisting of the non-singular
points of $X$ where $f$ is a submersion. Suppose
$Y$ is embedded in $X$, $f(Y)=0$.
To say that $f$ is $(t^0)$ regular over $Y$  means that $f$ satisfies the
$(w_f)$ condition for the pair $(X_f,Y)$ ($(w_f)$ for $(X,Y)$ means that
$\dist (T_0Y, \ker d(f\vert_X)) \leq C \dist (x, Y)$ for some constant
$C > 0$, cf. Definition 1.9).  If $X_f=X-Y$, then this condition ensures that
the family of functions on
$X$ parametrised by $Y$ induced from $f$ has a rugose trivialization (\cite{V},
\cite{T-W}).

We again have an integral closure version of this condition. The analogue
of 2.7 is:

\proclaim{Theorem 2.34} In the ambient $(t^r)$ set-up, suppose $X\subseteq
{\bk}^{n+k}$
is the germ of an analytic space at the origin,
$Y$ a submanifold embedded as
$\{0\}\times {\bk}^k$, $P$ an analytic direct
transversal to $Y$. For $r\ge 0$, $f$ is $(t^r)$
for the transversal $P$ iff
$$
\em_n^rJM_y(F)\subseteq \overline{
{\em_n}JM(F)_P+I(P)JM_y(F)}
$$
(take the integral closure inside $\Cal O_{X,0}^p$,
where $p$ is the number of components of $F$;
in the real analytic situation we take real
integral closure).\endproclaim

The analogue of 2.13 is similar; we use the same Grassmann modification of
$X$ and replace
$f$ by $\tilde f:=f\circ\beta$.

Remarkably, although ambient $(t^r)$ seems much stronger than asking that
$X$ and $X(f)$ are $(t^r)$, in the complex analytic case they are
equivalent. The precise statement in the two stratum case is:

\proclaim{Theorem 2.35} Suppose $X,0\subseteq {\bold {C}^{n+k}}$ is the germ
of a complex analytic space,
which contains $Y={\{0\}\times \bold{C}^k}$ as a stratum, $f:(X,Y)\rightarrow
(\bold {C},0)$, $X_f=X-Y$, and $P$ a direct transversal to $Y$. Then
$f$ satisfies $(t^r)$ for $P$ iff $X_0$ and $X(f)_0$ do.
\endproclaim

\demo{Proof} It is clear that if $f$ is ambient $(t^r)$ for $P$, then
$X_0$ and $X(f)_0$ (respectively the regular points of  $X$ and
$X(f)$) are $(t^r)$ for
$P$, because the limiting tangent hyperplanes
to the fibers of $f$ include the limiting tangent hyperplanes to $X$ and
to $X(f)$.

The other direction of the proof is by induction on $r$. Assume $X_0$ and
$X(f)_0$ are $(t^1)$ for $P$.  If we consider the Grassmann modification of
$X$, then $\tilde X_0$ and $\tilde X(\tilde f)_0$ are the complements of
$G_{n+k,n}$ in $\tilde X$ and $\tilde X(\tilde f)$ respectively. Hence, by the
proof of 2.13,
$\{\tilde X_0-\tilde X(\tilde f), \tilde X(\tilde f)_0, G_{n+k,n}\} $ form
a Whitney stratification of $\tilde X$ at $(j^1P,0)$. By the result of
Brian\c con, Maisonobe and Merle \cite{BMM}, this implies that
$(w_{\tilde f})$ holds at  $(j^1P,0)$ for $\{\tilde X_0, G_{n+k,n}\}$. Then
by the proof of 2.13, altered for
the ambient condition, $(t^1)$ holds for $f$ and $P$. Assuming $(t^r)$
holds for $X_0$ and $X(f)_0$ for $P$, again construct the Grassmann
modification, and apply 2.13 twice, using the
induction hypothesis.
\qed\enddemo

We will be using these results at the end of section 4.

\head \S 3.  Conditions $(a)$, $(t^1)$ and the aureole\endhead

Let $Y\subseteq X$  be analytic subsets
of $\bC^n$, $Y$ smooth, $S(X)\subseteq Y$. The relation
between the $(t^1)$ condition and
Whitney's  condition $(a)$ for $(X_0,Y)$
depends on the relation between the stratum
$Y$ and the {\it aureole} (\cite{H-L}, \cite{L-T}).  
Recall that the conormal modification
of $X$ is the closure in $X\times \bP^{n-1}$
of the set of hyperplanes containing
$T_xX$ at a regular point $x$.  Assume
$0\in X$. The fiber of the conormal
modification  over 0 will be denoted
$\C(X)$.  The aureole of $X$ at 0
is constructed in the complex analytic
case by blowing up the conormal
modification  of $X$ by the pull-back
of the maximal ideal of $X$
at 0, then projecting the irreducible
components of the exceptional
divisor $E$ of the blow-up to  the projectivised
tangent cone $\PTC(X)$ of $X$.  We will identify
subsets of $\PTC(X)$ with the
corresponding subcones of $\TC(X)$.
The images of these components consist
of a set of subcones; this set $\A$ is
the aureole. (The components of the
tangent cone are the largest members
of the aureole; any members which
are of lower dimension than the cone
are called {\it exceptional subvarieties}
of $\TC(X)$; exceptional subvarieties
cannot have dimension 0.)

As proved in \cite{L-T}, the aureole has the property that
the limiting tangent hyperplanes to
$X$ at 0 are exactly the planes tangent
to one of the elements of the aureole, i.e.
$\C(X)=\cup_{A\in\A}\, \C(A)$.
Thus the aureole tracks limiting tangent hyperplanes
and the limiting tangents of the curves
along which the limits are achieved. 
Without loss of generality, assume
$Y$ is linear. If $Y$ contains
an element of the aureole properly then Whitney $(a)$ will fail:
at a regular point of the contained
cone, the tangent space of the
contained cone will be a proper
subspace of $Y$, hence
there will be a tangent hyperplane $H$
tangent to the contained cone but not
containing $Y$; by the property of
aureoles mentioned above, $H\in\C(X)$,
hence $(a)$ fails. If there is an
$\ell\in \TC(X)$, $\ell\not\subseteq Y$,
a curve $\phi(t)$ in $X$ with $\phi(0)=0$,
$\phi(t)\in X_0$ for $t\ne 0$ and $\phi$
is tangent to $\ell$ at 0, and a family
of hyperplanes $H(t)\subseteq T_{\phi(t)}X_0$
such that $H=\lim H(t)$ does not contain $Y$,
than $(t^1)$ fails, and conversely
(see \cite{G4} Cor 2.8).  This allows
us to characterize both $(a)$ and
$(t^1)$ in terms of the aureole.

Let $\E$ denote the set of irreducible
components of the exceptional
divisor $E$.  Let $\pi_1$ denote the
projection from $E$ to the tangent cone;
so $\pi_1$ of the elements of $\E$
gives the members of $\A$; let $\pi_2$
denote the projection from $E$ to
$\C(X)$ (i.e. we ignore which line in
the tangent cone the hyperplanes come from).
Let $\E_1$  denote those elements
$C$ of $\E$ for which $\pi_2(C)\subseteq \C(Y)$;
let $\E_2=\E-\E_1$. 
\proclaim{Proposition 3.1} $X$ satisfies
$(a)$ over $Y$ at $0$ iff
$\E_2=\emptyset$.
\endproclaim
\demo{Proof} Note that $\C(X)=\cup_{C\in\E}\, \pi_2(C)$.
Thus $(a)$ holds iff $\C(X)\subseteq\C(Y)$
iff  $\pi_2(C)\subseteq\C(Y)$ for all $C\in\E$.
\enddemo

\proclaim{Proposition 3.2} $X$ satisfies
$(t^1)$ over $Y$ at $0$ iff
$\{C\in\E_2: \pi_1(C)\not\subseteq
Y\}=\emptyset$.
\endproclaim
\demo{Proof} By  Cor 2.8 of \cite{G4},
$(t^1)$ fails iff there is an $(\ell,H)\in E$
such that $\pi_1(\ell,H)=\ell\not\subseteq Y$
and $\pi_2(\ell,H)=H\notin \C(Y)$ iff
there exists $C\in\E_2$ such that
$\pi_1(C)\not\subseteq Y$.
\enddemo

The next result links Whitney $(a)$, the aureole
and $(t^1)$.

\proclaim{Theorem 3.3} Let $Y\subseteq X$ be analytic subsets of $\bC^n$, $Y$ smooth.
Assume $X$ satisfies $(t^1)$ over $Y$ at
$y\in Y$; then $X$ satisfies $(a)$
over $Y$ at $y\in Y$ iff there are no exceptional subvarieties of $\TC(X)_y$
properly contained in $T_yY$.
\endproclaim

\demo{Proof} Without loss of generality we
reduce to the case $y=0$ and $Y$ is linear.
Choose $C\in\E$. Since $(t^1)$
holds, either (i) $\pi_2(C)\subseteq \C(Y)$,
or (ii) $\pi_2(C)\not\subseteq
\C(Y)$ and
$\pi_1(C)\subseteq Y$. (Note that in (ii) the
inclusion must be proper; if not
$Y$ is an element of the aureole, and by
\cite {L-T}, $\pi_2(C)\subseteq \C(Y)$.)

Assume case (i); then
this $C$ presents no obstruction to $(a)$ holding.
Furthermore, if $A=\pi_1(C)$ is a proper analytic subset
of $Y$, then there will be  a tangent hyperplane
tangent to the $A$, but not containing $Y$, thus
contradicting the assumption.
Assume case (ii); then $(a)$ fails, and
$A=\pi_1(C)$ is a proper subvariety
of $Y$.
\qed
\enddemo

Observe that if $\dim Y=1$, no member of
the aureole can be a proper analytic subset
of $Y$, so $(t^1)$ implies $(a)$ (which was proved
in the subanalytic case by the second author (\cite{Tr1},
\cite{Tr4})).  If $\dim Y=2$,
Theorem 3.3 implies that $(t^1)$ is equivalent
to $(a)$ iff there are no exceptional lines
in $Y$.

We wish to give examples
for which $(t^1)$ holds but $(a)$ fails. Theorem 3.3 shows that it suffices
to find $X$ and $Y$ for which the aureole
$\A$ contains proper subvarieties
of $Y$, and for which all other members
$A$ of $\A$ have $\C(A)\subseteq \C(Y)$.

\example{Example 3.4} Let $I=(xyf(x,y,z,w))$,
$f=x^a+y^b+z^c+w^d$
with $2\le a<b<c<d$; let $X=V(I)$, $Y=S(X)$ the
$wz$ plane. The purpose of the $xy$ term is
to have $Y\subseteq X$. The aureole of $X$ is the union of the aureole of the three
components. The aureole consists of $x=0$, $y=0$
and the members of the aureole of $f=0$. It is
shown in \cite{G3} that the aureole of $f=0$
consists of a flag---the
linear spaces defined by
$x=0$,
$x=y=0$ and $x=y=z=0$; call these $A_1$,
$A_2$ and $A_3$, respectively. Then $A_3$ is
a proper algebraic subset of $Y$, and
$A_1$ and $A_2$ satisfy
$\C(A)\subseteq \C(Y)$.
Thus  $(t^1)$ holds but $(a)$ fails.
\par
Without going into the details, we consider what
happens in this example if some of the strict inequalities
are replaced with equalities.  If $2\le a=b<c<d$,
$\A(f=0)=\{\{x^a+y^a=0\},\{x=y=0\},\{x=y=z=0\}\}$, where
$x^a+y^a=0$ is the union of $a$ hyperplanes
intersecting in $x=y=0$; $(t^1)$ holds but $(a)$ fails.
If $2\le a<b=c<d$, $\A(f=0)=\{\{x=0\},
\{x=0,y^b+z^b=0\},\{x=y=z=0\}\}$,
where
$x=0,y^b+z^b=0$ is the union of $b$ planes
intersecting in $x=y=z=0$; $(t^1)$ fails.
If $2\le a=b=c<d$,  $\A(f=0)=\{
\{x^a+y^a+z^a=0\},\{x=y=z=0\}\}$; $(t^1)$ fails.
\endexample

In the real case, the structure of the
conormal modification is much more complicated
(see \cite{O-W} for the case of surfaces
in $\bR^3$).  However we state one collection
of examples whose properties are not too hard
to check.

\example{Example 3.5} Now we work in $\bR^4$.
Let $I=(xyf(x,y,z,w))$, $f=x^a+y^b-z^c+w^d$
with $2\le a= b<c<d$, $a,b,c,d$ even;
let $X=V(I)$, $Y$ the $wz$ plane.  The geometric
tangent cone of $f=0$ is $Y$, so direct
transversals to $Y$ intersect $f=0$ only at 0.
$\C(X)$ consists of all hyperplanes containing
$x=y=z=0$.
Thus  $(t^1)$ holds but $(a)$ fails.
In \cite{Tr4} the second author previously gave semialgebraic examples in
$\bR^4$ showing that $(t^1)$ can hold without $(a)$.
\endexample 
\head \S 4. The Principle of Specialization
of Integral Dependence and
the $(t^r)$ conditions\endhead

So far we have not used very strongly the
algebraic character of the theory of
integral closure.  However, if we restrict
ourselves to complex analytic sets, there is a
remarkable result which comes from this side of the theory.
In many applications we have a
family of sets, and a sheaf of modules on each set.
These modules often come from
specializing a sheaf of modules contained
in a locally  free sheaf
on the total space of the family to the
members of the family. Often we know
that some element of the free sheaf
is in the integral closure of the modules
on the members of the family generically. Is it in the integral closure of the
member-module for all members? Is it in the
integral closure of the sheaf of modules which is
defined on the total space? This type of question often
arises in mathematics, and the answer is usually no.

The theory of integral closure in some
cases associates an invariant to the
member-modules called their multiplicity.
If these multiplicities
are constant then the answer to both questions
is yes! This is called the
Principle of Specialization of Integral Dependence,
first stated for ideals by Bernard Teissier
(see \cite{T1}, 3.2, p.330 and \cite{T2}, App. I).

For all statements and results of this section,
we will assume $\bk=\bC$.
If one considers the maximal ideal $\em_x$
of a local ring  ${\Cal O}_{X,x}$,
where $X$ is equidimensional at $x$, then
the multiplicity of $\em_x$ is
just the multiplicity of $X$ at $x$---the
number of points near $x$ in which a
generic linear space of dimension complementary
to $X$ will intersect $X$.  If one considers an
ideal $I$ of finite colength
(the colength of a module $M\subseteq
{{\Cal O}_X}^p$ is the
dimension of the quotient module
$\displaystyle{{{\Cal O}_X}^p}\over M$ as a
complex vector space), one can pick
a set of generators
of $I$, and use these to define a map-germ
$f$ with image $f(X)$.  The
multiplicity of $I$  is just the product
of the multiplicity of $f(X)$ and the
degree of $f$.

The multiplicity of a module is defined
for all modules of finite colength
(\cite{B-R} and \cite{K-Th}).
If we are working over a ring which
is Cohen-Macaulay then the multiplicity $e(M;X)$
of  $M\subseteq {\Cal O}^p_X$ can be obtained
by taking the colength of the submodule
generated by $d+p-1$ generic elements of
$M$, where $d$ is the
dimension of $X$ (\cite{G2}, proposition 2.3,
and the fact that any submodule generated by these
generic elements has the same
multiplicity as $M$ \cite{K-Th}).
Before stating the Principle of
Specialization of Integral Dependence precisely,
we give an example to make these ideas concrete.

\example{Example 4.1} Consider the ideal
$I=(x^2,y^2,txy)$ in ${\Cal O}_3$; we
think of ${\bC}^3$ as the total space
of a family of planes parametrised by
$\bC$ embedded as the $t$-axis. The ideal $I$ induces
a family of ideals of finite colength $I_t$. Clearly the term $xy\in I_t$ for
$t\ne 0$.  The multiplicity
of $I_t$ is 4 for all values of $t$, so we know that
$xy\in{\overline I}$.
\endexample

If we have a mapping between two sets
$X$ and $Y$, and the restriction of
a sheaf of modules to the fiber
over $Y$  is supported at a finite number of points,
then we denote by $e(y)$ the
sum over the points
in the fiber of y of the multiplicities
of the stalks of the sheaf.

\proclaim {Theorem 4.2}
(Specialization of integral dependence) \newline
Assume that $X$ is equidimensional, and that
$F\:(X,x_0)\to(Y,y_0)$ is a map
of germs of analytic spaces.
Assume that the fibers $X(y)$ are
equidimensional of the same dimension
$d$, $d$ at least 1, and that $Y$ has
dimension at least 1. Let $M$
be a coherent submodule of ${\Cal O}^p_X$.  Set
$S:=Supp({\displaystyle{\Cal O_X^p}\over M})$,
and assume
that $S$ is finite over $Y$. Assume $y\mapsto e(y)$
is constant on $Y$.  Let $h$ be a section of
$\Cal O^p_X$ whose image in
$\Cal O^p_{X_y}$ is integrally
dependent on the image of $M$ for $y$ in
a Zariski open set.
Then (after $X$ and $Y$ are replaced
by neighborhoods of $x_0$ and
$y_0$ if necessary) $h$ is integrally
dependent on $M$.\endproclaim
\demo {Proof} The proof appears in \cite{G-K}.
\qed\enddemo
Now suppose $X=F^{-1}(0)$; we use the multiplicity of the
module $JM(F)_{P_r}$ ($P_r$ is defined in (2.2))
thought of as a submodule
of a free $\Cal O_{X\cap P_r,0}$ module to control the
Whitney equisingularity of a family of
generic intersections with $X$. It is not surprising that
the family of generic intersections
is Whitney equisingular; what is surprising is that only
one invariant is needed. In general,
$d+1$ invariants of the Jacobian module are needed where $
d$ is the dimension of the elements of the family (\cite{G2}).

\proclaim{Theorem 4.3} (Genericity Theorem)
Suppose $X\subseteq \bC^n\times \bC^k$ is a
$d$-dimensional complete intersection, and
suppose $P$ is a direct transversal to
$Y=\{0\}\times \bC^k$, $P\cap S(X)=\{0\}$.
Then $e(JM(F)_{P_r}; P_r\cap X)$ is the
generic (finite) value of this invariant
for the $r$-jet of $P$ among all transversals
with the same $(r-1)$-jet iff $(t^r)$ holds for $P$.
\endproclaim
\demo{Proof} The proof is by induction on $r$. Note that $(t^r)$ holds for $P$ iff $(t^r)$
holds for $P_r$, so we can reduce to the $
r$-jet of $P$. Further, since $P\cap S(X)=\{0\}$,
it follows that the dimension of $S(X)$
is less than or equal to that of $Y$, hence the
generic $r$-jet intersects $S(X)$ only at $0$.
(This claim is obvious for the generic $1$-jet;
if we fix the $r-1$ jet and apply the Grassmann
modification $r-1$ times, varying the
$r$-jet amounts to varying the family of planes
in the $r-1$ modification.)

Let $r=1$. We can view the modification $\tilde X$
as the family of intersections $\{X\cap P\}$
parameterised by the Grassmannian where
$P$ is a plane which is a direct transversal.
In the proof of Theorem 4.3 of \cite{G4}, it
is shown that the planes for which
the Jacobian module of $\{X\cap P\}$ has the
generic value of the multiplicity are the                                                                      same as
the planes in which $JM(F)_{P}$ has the generic value.
Call the Z-open subset of the Grassmannian which
parameterises these generic planes $U$. Since the
multiplicity of the Jacobian module is
constant along $U$, it follows that there exists
some open neighborhood $W'$ of $U$ in $\tilde X$
which consists only of points in $U$ and in
$\beta^{-1}(X-S)$. This follows because there is a p
erhaps smaller subset of $U$ which consists of planes intersecting
the tangent cone of $S$ only at the origin. For these
planes there are no points of $S$ close to zero; but
since the multiplicity of the Jacobian module is constant
along $U$, there can be no points of $S$ for any of these
planes as these would cause the multiplicity to jump. This
shows that the planes parametrised by $U$ satisfy the part of the
$(t^1)$ condition concerned with $S$. The proof that the
part of the $(t^1)$ condition concerning $X-S$ holds
follows from the  proof of Theorem 4.3
of \cite{G4}. If we assume that the $(t^1)$ holds for
$P_1$, then the proof that the value of  $JM(F)_{P}$
is the generic value follows as in \cite{G4}.

The idea of the proof is to show
that the corresponding rings
and modules at the level of $X\cap P_r$
and at the level of
$\tilde X \cap \tilde P_r$ are isomorphic. This will show that the multiplicities
of the corresponding modules are the same. The
result for $(t^{r-1})$ at the level of
$\tilde X$ will then imply the result
for $(t^r)$ at the level of $X$. The
hypotheses on $P$ and $S(X)$ are to
ensure that $e(JM(F)_{P_r}; P_r\cap X)$
is generically finite.  Suppose $P_r$
is given by $\{y_i-p_i(x)=0\}$, $p_i$ a
polynomial of degree $r$.  Suppose
$p_i(x)=\sum p_{i,j}(x) x_j$, giving a lift
$\tilde P$ of $P_r$ with equations
$\{a_{i,j}-p_{i,j}(x)=0\}$.

We can define maps $p=(p_1,\dots,p_k)$
and $\tilde p=(p_{i,j})$ such that $P_r$
and $\tilde P$ are isomorphic to $\bC^n$
by the graph maps $p_\Gamma$ and
$\tilde p_\Gamma$, and we have that
$\beta\circ \tilde p_\Gamma=p_\Gamma$. Then the rings $\Cal O_{\tilde X\cap
\tilde P,P_1}$ and $\Cal O_{X\cap P_r,0}$
are isomorphic, since
$$
\tilde p^*((a_{i,j}-p_{i,j}), \beta^*(I(X)))=
\tilde p^*\beta^*(I(X))=p^*I(X)
=p^*(I(X),(y_i-p_i)).
$$

Now we consider the modules.  If we
restrict the module used in the
formulation of the $(t^r)$ condition
to $X\cap P_r$ we see that it is generated by
$\em_n$ times
$$
\left\{ \pd F{x_j} +\sum_{i=1}^k \pd
{p_i}{x_j} \pd F{y_i} : 1\le j\le n\right\},
$$
where $X=F^{-1}(0)$.

The corresponding module which controls
the $(t^{r-1})$ condition on $\tilde X$
when restricted to $\tilde X\cap \tilde P$
is generated by  $\em_n$ times
$$
\left\{ \pd G{x_j} +\sum_{i,l}
\pd {p_{i,l}}{x_j} \pd
G{a_{i,l}} : 1\le j\le n\right\},
$$
where $G=F\circ\beta$.

Now, $\pd G{a_{i,l}}=x_l \pd F{y_i}
\circ\beta$, while $\pd G{x_j}=\pd
F{x_j}\circ\beta + \sum_i a_{i,j}\pd
F{y_i}\circ\beta$.  So, the generators of
the module on $\tilde X$ become $\em_n$ times
$$
\left\{\pd F{x_j}\circ\beta +\sum_i a_{i,j}
\pd F{y_i}\circ\beta +\sum_{i,l} \pd
{p_{i,l}}{x_j} x_l \pd F{y_i}\circ\beta\right\}.
$$
Now, $\pd{p_i}{x_j}=\pd{\sum p_{i,l}x_l}{x_j} =
(\sum_i \pd{p_{i,l}}{x_j} x_l)+p_{i,j}$ by
summing over $l$, so we can rewrite the
generators of the module on $\tilde X$
restricted to $\tilde X\cap \tilde P$ as
$\em_n$ times
$$
\multline
\left\{\pd F{x_j}\circ\beta +\sum_i p_{i,j}
\pd F{y_i}\circ\beta + \sum_{i=1}^k \pd
{p_i}{x_j}\pd F{y_i}\circ\beta- \sum_i p_{i,j}
\pd F{y_i}\circ\beta
\right\}=\\
\left\{\pd F{x_j}\circ\beta +\sum_{i=1}^k
\pd {p_i}{x_j}
\pd F{y_i} \circ\beta \right\}.
\endmultline
$$

Since these generators pull back to
the same elements when composed with
$p_n$ and $\tilde p_n$, the corresponding
quotient modules are isomorphic, hence
have the same multiplicity.
Varying the $r$-jet of $P$ with $(r-1)$-jet
fixed amounts to varying the $(r-1)$-jet
of $\tilde P$ with the $(r-2)$-jet fixed.
So by induction we are done.
\qed\enddemo

\proclaim{Corollary 4.4} If we fix an
$(r-1)$-jet, then the set of $r$-jets which
satisfy $(t^r)$ is Zariski open and
dense.
\endproclaim

\proclaim{Corollary 4.5} Consider
the family obtained by fixing an $(r-1)$-jet and
varying the $r$-jet. Then the family has
Verdier equisingular intersection with
$X$ exactly for those parameter
values for which the multiplicity takes its
generic value.
\endproclaim

\demo{Proof}
\comment
The proof
is similar to the proof of the Genericity
Theorem.The case where $r=1$ follows from the
proof of (4.3) of \cite{G4}, and the general case
follows from the proof of the Genericity
Theorem, because we can lift all of our hypotheses
to the Grassmann modification.
\endcomment
This follows from Theorem 4.3 and Corollary 2.28.\qed
\enddemo

By Theorem 4.3 we know that if $P$ is a direct
transversal such that the multiplicity,
$e(JM(F)_{P_r},\Cal O_{X\cap P_r})$ has the
generic value $e$ among all such modules
associated to $r$-jets with the same $(r-1)$-jet as
$P$, then $X$ is $(t^r)$ for $P_r$, hence for $P$.
Call the multiplicity $e(JM(F)_{P},\Cal
O_{X\cap P})$, the associated multiplicity of the
pair $(X,P)$. It is not difficult to see that if
$X$ is $(t^r)$ for $P$, then the associated
multiplicity is the generic value $e$.   For the
deformation $P(t)$ to $P_r$ has constant
$r$-jet, and $X$ is $(t^r)$ for $P_r$, hence
by Corollary 2.28 the deformation is a
Verdier equisingular family of isolated complete
intersection singularities,
hence by Theorem 1 and Proposition 2.6  of \cite{G2},
the multiplicity of the Jacobian module of $X\cap
P(t)$ is constant. The multiplicity of this
Jacobian module is the same as the associated
multiplicity, because the quotients of these two
modules are isomorphic as ${\Cal O}_{X\cap
P}$-modules. If $X$ is not $(t^r)$
for $P$, the associated
multiplicity of $(X,P_r)$ must be greater than
the generic value by Theorem 4.3. It may be that
the associated multiplicity of $(X,P)$ is less
than that of $(X,P_r)$; can it be the generic
value of $e$ or even less? The next lemma deals
with this question.

\proclaim{Lemma 4.6} Suppose
$(t^r)$ fails for a direct transversal $P$.
Then the value of $e(JM(F)_P,\allowmathbreak\Cal
O_{X\cap P})$ is strictly greater than $e$, the
generic value among all $P'_r$ such that
$j^{r-1}P=j^{r-1}P'$.
\endproclaim

\demo{Proof}  By
the modification theorem (2.13) and the techniques used
in the Genericity Theorem (4.3) it suffices to prove the
lemma for $(t^1)$.
Let the equations
for $P$ be $\{y_i-\sum p_{i,j}(x) x_j=0\}$. Consider
the family of direct transversals given by
$\{y_i-\sum (a_{i,j}+p_{i,j}(x))x_j=0\}$ in
$\bC^{n+k}$.  Denote by $P^a$ the
transversal determined
by $a=(a_{i,j})$. We can study the resulting
family of sections by the following
construction. Let $\MAP{\tilde\beta}{\bC^n\times
\bC^{nk}}{\bC^{n+k}}$ be given by
$$
\tilde\beta(x_1,\dots,x_n,a_{1,1},
\dots,a_{k,n})=(x_1,\dots,x_n,\mathop{\sum}\limits_i
(a_{1,i}+p_{1,i}(x))x_i,\dots,\mathop{\sum}\limits_i
(a_{k,i}+p_{1,i}(x)) x_i).
$$
Let $\tilde X$ denote
the set defined by $F\circ\tilde\beta$. Given
$\phi :{\bC},0 \rightarrow \tilde X,P\times
\{0\}$, then $\tilde\beta\circ\phi$ is tangent to
$P$ at the origin, and only curves tangent to $P$
at zero lift to $\tilde X$. Tracing through the
proof of Theorem 2.13, it is easy to see that
$(t^1)$ holds for $P^a$ iff
$\tilde\beta^*JM(F)_{P^a}$ is a reduction of
$\tilde\beta^*JM(F)$. Since $(t^1)$ holds for $P^a$
for generic $a$, this inclusion is true
generically. At these generic $a$ values, the
value of the associated multiplicity of $(X,P^a)$
is the generic value for 1-jets at the origin of
$X$, since $(t^1)$ holds. This shows that the
associated multiplicity of $(X,P)$ must be at least
the generic value since by \cite{G-K},
Proposition 1.1, the multiplicity is upper
semicontinuous.  Now suppose the associated
multiplicity has the generic value of the
associated multiplicities of 1-jets.  Then on
$\tilde X$ in a neighborhood of $(0,0)$ the
associated multiplicities are constant, and the
Principle of Specialization of Integral Dependence
shows that $\tilde\beta^*JM(F)_{P^0}$ is a
reduction of  $\tilde\beta^*JM(F)$. Thus $P$ is
$(t^1)$, which is a
contradiction.\qed
\enddemo

\proclaim{Theorem 4.7}
Suppose we have a family $P(t)$ of direct
transversals with fixed $(r-1)$-jet, and for some
parameter value $e(JM(F)_{P(t)},\Cal O_{X\cap
P(t)})$ has the generic value $e$ of
the multiplicity.  Then the family has Verdier
equisingular intersection with $X$
on the Zariski open, dense subset of
values where $e(JM(F)_{P(t)},\Cal O_{X\cap
P(t)})=e$ and it fails to have Verdier equisingular
intersection with $X$ where $e(JM(F)_{P(t)},\Cal
O_{X\cap P(t)})>e$.
\endproclaim

\demo{Proof}  From
the hypotheses we know that the set of parameter
values where the associated multiplicity has the
generic value $e$ is a
non-empty Zariski open set $U$. By
Lemma 4.6, we know that $(t^r)$ holds holds for
$P(t)$ exactly on this Zariski open set. Then the
restriction of the family of sections to $U$ is
Verdier equisingular.  At points on the complement
of $U$, the associated multiplicity jumps, so the
multiplicity of the Jacobian module of the
section jumps, and Verdier equisingularity fails
by Theorem 1 of \cite{G2}.\qed
\enddemo

We have seen that if we
fix $P_{r-1}$ and vary the $r$ jet, then only one
invariant, the associated multiplicity, need be
constant to ensure Verdier equisingularity,
provided this invariant takes its generic value.
What if the value of the associated
multiplicity is greater than the generic value?
If we have some control over the plane sections
of $P$, then it is still possible to say
something. Given an isolated complete
intersection singularity $X,0\subset\bC^n$ we
can consider the $\mu_*$ sequence of $X$. This is
the sequence of Milnor
numbers $\mu_0(X)=\mu(X)$,\dots,
$\mu_i(X)=\mu(X\cap H_i)$,\dots,
$\mu_{d+1}(X)=1$, where $d=\dim X$
and $H_i$ is a generic linear
space of codimension $i$. Sometimes authors index
by the dimension of $H$, so $\mu^d(X)=\mu(X\cap
H^d)$. There is a nice connection between the
associated multiplicity of $(X,P)$ and the $\mu_*$
numbers of $X\cap P $. If  $X\cap P $ is an isolated complete
intersection singularity,
then the associated multiplicity of $(X,P)$ is the
sum of $\mu(X\cap P)$ and $\mu ((X\cap
H)\cap(P\cap H))$. This follows from the fact that
the associated multiplicity is the multiplicity
of the Jacobian module of $X\cap P$, Proposition
2.6 of \cite{G2}, and the theorem of L\^e
and Greuel (cf. p. 211 of \cite{G2}, material after
1.2). Coupling this connection with Theorems 1 and
2 of \cite{G-2} shows that a family of isolated complete
intersection singularities is
Verdier equisingular iff the $\mu_*$ invariants
are independent of parameter.

\proclaim{Theorem
4.8} Suppose $X$ is a complete intersection
defined by $F$. Suppose $P(t)$ is a family of
direct transversals which intersect $X$ in the
expected dimension, $P(t)\cap S(X)=0$, and the
associated multiplicity of $(X,P(t))$ is constant.
Suppose further that $\mu_1(X\cap
P(t))=\mu_{k+1}(X)$, and $\mu_2(X\cap
P(t))=\mu_{k+2}(X)$.  Then the family $\{X\cap
P(t)\}$ is
Verdier equisingular.
\endproclaim

\demo{Proof}
Because $\mu_1(X\cap P(t))=\mu_{k+1}(X)$,
and $\mu_2(X\cap P(t))=\mu_{k+2}(X)$, it follows
that the associated multiplicity of $(X\cap H,
P(t)\cap H)$ is the generic value for planes of the
same dimension as $(P(t)\cap H)_1$, $H$ a generic
hyperplane. It follows from Lemma 4.6, that
$X\cap H$ is $(t^1)$ for $P(t)\cap H$, hence $(X\cap
H)\cap(P(t)\cap H)$ and its sections are Verdier
equisingular. Since the associated multiplicity
of $(X,P(t))$ is constant, it follows that the Milnor
number of $X\cap P(t)$ is constant as well.
Since the $\mu_*$ invariants are constant, the
family $X\cap P(t)$ is Verdier equisingular as
well.\qed
\enddemo

We can also prove analogues of the theorems of this section for the
ambient $(t^r)$ conditions and the  $(w_f)$ condition. In the ambient
$(t^r)$ set-up,
we use the multiplicity of  $e(JM(F)_{P_r},\Cal O_{X\cap P_r})$. If
$X=\bC^n$, then this is just the Milnor number of $f|{P_r}$. A typical
result is the following:

\proclaim{Theorem 4.9} Suppose we have a family $P(t)$ of direct
transversals with fixed $(r-1)$-jet, and for some
parameter value $\mu(f|P(t))$ has the generic value $\mu$ of
the Milnor number.  Then the family of functions has a rugose
trivialization
on the Zariski open, dense subset of
values where $\mu(f|P(t))=\mu$ and it fails to have such a trivialization
where $\mu(f|P(t))>\mu$.

\endproclaim
 Again in the hypersurface case, we
can do better than 4.8  because we are able to operate
in the ambient space, not just on the variety. Again the invariant of
interest is the
multiplicity of the ideal $J(f|P)$which is just the Milnor number of
$f|P$.

\proclaim{Theorem 4.10} Suppose $X$ is a
hypersurface defined by
$f$, $X\subset\bC^n\times\bC^k$. Suppose $P(t)$ is
a family of direct transversals which intersect
$X$ in the expected dimension, and $P(t)\cap
S(X)=0$. Suppose also $\mu(X\cap P(t))$ is
constant and $\mu_1(X\cap P(t))=\mu_{k+1}(X)$. Then
the family $\{X\cap P(t)\}$ is Verdier
equisingular.
\endproclaim

\demo{Proof} Let $H_t$
be a generic hyperplane so that $\mu(X\cap
P(t)\cap H_t)=\mu_{k+1}(X)$. Note that, for a
generic plane $H'$ of codimension
$k+1$, $\mu_{k+1}(X)=\mu(X\cap H')$ iff
$\overline{J(f)}|_{H'} = \overline{J(f)_{H'}}$,
and in the above case $e(J(f)|_{H'}) =
\op{col}(J(f)_{H'},I(H'))$ ($\op{col}$ is the colength).
Now $\mu(X\cap P(t)\cap
H_t)=\op{col}(J(f)_{P(t)\cap H_t},I(P(t)\cap
H_t))$. Here we use the ambient genericity lemma:
if $\mu_{k+1}(X)=\mu(X\cap H')$, then the ambient
$(t^1)$ holds for $H'$, i.e.
$\em_nJ_y(F)\subseteq\overline{\em_nJ(F)_{H'}+I(H')J_y(F)}$.
Then
a proof similar to that of Lemma 4.6 shows that
$\mu(X\cap P_1(t)\cap H_t)=\mu_{k+1}(X)$. At this
point the proof proceeds as in Theorem
4.8.\qed
\enddemo

\example {Example 4.11} Suppose
$X$ is a hypersurface in $\bC^n$. If $n=3$ and
$P_t=H_t$ is a family of hyperplanes with
$H_t\cap X$ having constant Milnor number, then the
family is Verdier equisingular because it is a
family of plane curves. If $n=4$, then we need
$\mu(H_t\cap X)$ and $\mu_1(H_t\cap X)$
constant. If $n=5$, then we need $\mu(H_t\cap X)$
constant and
$\mu_1(H_t\cap X)=\mu_2(X)$.
\endexample

\address
Department of Mathematics, Northeastern University, Boston, Mass., U.S.A.
\endaddress
\address
LATP (UMR 6632), Centre de Math\'ematique et Informatique, Universit\'e
de Provence, 39 rue Joliot-Curie, 13453 Marseille, France
\endaddress
\address
Department of
Mathematics, University of Hawaii at Manoa, 2565 The Mall, 96822 Honolulu,
Hawaii, U.S.A.
\endaddress
\enddocument